\documentclass[12pt,leqno]{article}
\usepackage{amsfonts}
\usepackage{amsmath, amsthm, amsfonts, amssymb, color}
\usepackage{mathrsfs}
\usepackage{color}
\usepackage[notcite,notref]{showkeys}

\theoremstyle{definition}

\newcommand{\scr}[1]{\mathscr #1}
\definecolor{wco}{rgb}{0.5,0.2,0.3}

\numberwithin{equation}{section} \theoremstyle{remark}

\newcommand{\ua}{\uparrow}

\title{{\bf Functional Inequalities for Weighted Gamma Distribution on the Space of  Finite  Measures}\footnote{Supported in
 part by  NNSFC (11771326, 11831014).} }
\author{{\bf   Feng-Yu Wang$^{a,b)}$  }\\
\footnotesize{$^{a)}$ Center for Applied Mathematics, Tianjin University, Tianjin 300072, China}\\
 \footnotesize{$^{b)}$ Department of Mathematics,
Swansea University, Bay Campus, Swansea SA1 8EN, United Kingdom}\\
\footnotesize{    wangfy@tju.edu.cn, F.-Y.Wang@swansea.ac.uk}}
\begin{document}
\allowdisplaybreaks
\def\R{\mathbb R}  \def\ff{\frac}   \def\B{\mathbf B}
\def\N{\mathbb N} \def\kk{\kappa} \def\m{{\bf m}}
\def\ee{\varepsilon}\def\ddd{D^*}
\def\dd{\delta} \def\DD{\Delta} \def\vv{\varepsilon} \def\rr{\rho}
\def\<{\langle} \def\>{\rangle} \def\GG{\Gamma} \def\gg{\gamma}
  \def\nn{\nabla} \def\pp{\partial} \def\E{\mathbb E}
\def\d{\text{\rm{d}}} \def\bb{\beta} \def\aa{\alpha} \def\D{\scr D}
  \def\si{\sigma} \def\ess{\text{\rm{ess}}}\def\s{{\bf s}}
\def\beg{\begin} \def\beq{\begin{equation}}  \def\F{\scr F}
\def\Ric{\mathcal Ric} \def\Hess{\text{\rm{Hess}}}
\def\e{\text{\rm{e}}} \def\ua{\underline a} \def\OO{\Omega}  \def\oo{\omega}
 \def\tt{\tilde}\def\[{\lfloor} \def\]{\rfloor}
\def\cut{\text{\rm{cut}}} \def\P{\mathbb P} \def\ifn{I_n(f^{\bigotimes n})}
\def\C{\scr C}      \def\aaa{\mathbf{r}}     \def\r{r}
\def\gap{\text{\rm{gap}}} \def\prr{\pi_{{\bf m},\varrho}}  \def\r{\mathbf r}
\def\Z{\mathbb Z} \def\vrr{\varrho} \def\ll{\lambda}
\def\L{\scr L}\def\Tt{\tt} \def\TT{\tt}\def\II{\mathbb I}
\def\i{{\rm in}}\def\Sect{{\rm Sect}}  \def\H{\mathbb H}
\def\M{\scr M}\def\Q{\mathbb Q} \def\texto{\text{o}} \def\LL{\Lambda}
\def\Rank{{\rm Rank}} \def\B{\scr B} \def\i{{\rm i}} \def\HR{\hat{\R}^d}
\def\to{\rightarrow}\def\l{\ell}\def\iint{\int}
\def\EE{\scr E}\def\Cut{{\rm Cut}}\def\W{\mathbb W}
\def\A{\scr A} \def\Lip{{\rm Lip}}\def\S{\mathbb S}
\def\BB{\scr B}\def\Ent{{\rm Ent}} \def\i{{\rm i}}\def\itparallel{{\it\parallel}}
\def\g{{\mathbf g}}\def\Sect{{\mathcal Sec}}\def\T{\mathcal T}\def\BB{{\bf B}}
\def\f{\mathbf f} \def\g{\mathbf g}\def\BL{{\bf L}}\def\MM{\mathbb M} \def\BG{{\scr G}}
\def\Bd{{\nn^{ext}}} \def\BdP{\nn^{ext}_\phi} \def\Bdd{{\bf \dd}} \def\Bs{{\bf s}}\def\BD{{\bf D}}\def\GA{\scr A}
\def\Bg{{\bf g}}  \def\Bdd{{\bf d}} \def\supp{{\rm supp}}
\def\ddiv{{\rm div}}\def\osc{{\bf osc}}\def\1{{\bf 1}}
\def\BGG{{\bf\Gamma}}
\maketitle

\begin{abstract}  Let $\MM$ be the space of finite measures on a locally compact Polish space, and let  $\BG$ be   the  Gamma distribution on $\MM$  with intensity  measure $\nu\in \MM$. Let $\nn^{ext}$ be the extrinsic  derivative  with tangent bundle $T\MM= \cup_{\eta\in\MM}   L^2(\eta)$, and let $\GA: T\MM\to T\MM$ be measurable such that $\GA_\eta$ is a positive definite linear operator on $L^2(\eta)$ for every $\eta\in \MM$. Moreover, for a measurable function  $V$ on $\MM$, let $\d\BG^V= \e^V\d\BG$.  We investigate the Poincar\'e,  weak Poincar\'e and super Poincar\'e inequalities   for the Dirichlet form $$\EE_{\GA,V}(F,G):= \int_\MM \<\GA_\eta\nn^{ext}F(\eta), \nn^{ext}G(\eta)\>_{L^2(\eta)}\, \d\BG^V(\eta),$$ which characterize various properties of the associated Markov semigroup.   The main results are   extended to the space of finite signed measures.\end{abstract} \noindent
 AMS subject Classification:\  60G57, 60G45, 60H99.   \\
\noindent
 Keywords:  Extrinsic derivative,    weighted Gamma distribution, Poincar\'e inequality, weak Poincar\'e  inequality, super Poincar\'e inequality.
 \vskip 2cm

\section{Introduction}
Let $\MM$ be the class of finite measures on a locally compact Polish space $E$, which is again a Polish space under the weak topology. Recall that
a sequence of finite measures $\eta_n\to\eta$ weakly if $\eta_n(f)\to\eta(f)$ for $f\in C_b(E)$, where and in what follows, for a measure $\eta$ we denote
\beq\label{*00} \eta(f):= \int f\d\eta,\ \ f\in L^1(\eta).\end{equation}
Since $M$ is locally compact, the Borel $\si$-algebra $\B(\MM)$ induced by the weak topology coincides with that induced by the vague topology.
Let $\nu\in \MM$ with $\nu(E)>0$.  The Gamma distribution   $\BG$ with intensity measure $\nu$ is the unique probability measure on $\MM$ such that for any finitely many disjoint measurable subsets $\{A_1,\cdots, A_n\}$ of $E$,
$\{\eta(A_i)\}_{1\le i\le n}$ are independent Gamma random variables with shape parameters $\{\nu(A_i)\}_{1\le i\le n}$ and scale parameter $1$; that is,
\beq\label{*0}  \int_{\MM} f(\eta(A_1),\cdots, \eta(A_n))\BG(\d\eta)
 = \int_{[0,\infty)^n} f(x_1,\cdots, x_n)\prod_{i=1}^n \gg_{\nu(A_i)}(\d x_i),\ \ f\in \B_b(E), \end{equation}  where $\B_b(E)$ is the class of bounded measurable functions on $E$, for a constant $r>0$
 \beq\label{*0'}  \gg_r(\d s):= 1_{[0,\infty)}(s) \ff{s^{r-1}\e^{-s}}{\GG(r)} \, \d s,\ \  \GG(r):=\int_0^\infty x^{r-1}\e^{-x}\d x,\end{equation} and   $\gg_0:= \dd_0$ is the Dirac measure at point $0$.
It is well known that $\BG$ is concentrated  on the class of  finite discrete measures
 $$\MM_{dis}:=\Big\{\sum_{i=1}^\infty s_i\dd_{x_i}:\ s_i\ge 0, x_i\in E, \sum_{i=1}^\infty s_i<\infty\Big\}.$$

Consider the weighted Gamma distribution $\BG^V(\d\eta):= \e^{V(\eta)}\BG(\d\eta)$, where $V$ is a measurable function on $\MM$. We will investigate functional inequalities for the Dirichlet form  induced by  $\BG^V(\d\eta)$  and a positive definite linear map $\scr A$ on the tangent space of the  extrinsic derivative. See \cite{HKLV} and references therein for Dirichlet forms induced by both extrinsic and intrinsic derivatives, where the intensity measure $\nu$ is the Lebesgue measure on $\R^d$ such that the Gamma distribution $\BG$ is concentrated on the space of infinite Radon measures on $\R^d$. In this paper, we only consider finite intensity measure $\nu$.


 \beg{defn}[\cite{ORS}]   A measurable real function $F$ on $\MM$ is called extrinsically differentiable at $\eta\in \MM$,    if
 $$\nn^{ext} F(\eta)(x):= \ff{\d}{\d s} F(\eta+s\dd_x)\Big|_{s=0}\ \text{exists\ for\ all} \   x\in E,$$ such that
 $$\|\nn^{ext} F(\eta)\|:= \|\nn^{ext} F(\eta)(\cdot)\|_{L^2(\eta)} <\infty.$$    If $F$ is   extrinsically differentiable at all $\eta\in \MM$,  we denote $F\in \D(\nn^{ext})$ and call it extrinsically differentiable on $\MM$.
 \end{defn}

Regarding $L^2(\eta)$ as the extrinsic tangent space at $\eta\in\MM$,   we define the directional derivatives   by   $$\nn^{ext}_\phi F(\eta):= \<\nn^{ext}F(\eta), \phi\>_{L^2(\eta)}=\eta\big(\phi \nn^{ext}F(\eta)\big),\ \ \phi\in L^2(\eta).$$
When $\phi$ is bounded, this coincides with the directional derivative under multiplicative actions:
    $$\nn^{ext}_\phi F(\eta)=\ff{\d}{\d s} F(\e^{s \phi}\eta)\Big|_{s=0} = \ff{\d}{\d s} F((1+s\phi)\eta)\Big|_{s=0},\ \ \phi\in \B_b(E),$$
where $h\eta$ for $h\in \B_b(E)$ is a finite signed measure given by
$$(h\eta)(A):= \eta(1_Ah)=\int_A h\,\d\eta,\ \ A\in\B(E).$$

To introduce the Dirichlet form induced by the extrinsic derivative and the weighted Gamma distribution $\BG^V$, we consider the class  $\F C_0^\infty$, which consists of   cylindrical functions functions of type
$$  F(\eta):=f(\eta(A_1),\cdots, \eta(A_n)),\ \ n\ge 1, f\in C_0^\infty(\R^n),  \{A_i\}_{1\le i\le n}\in\scr I(E),$$ where $\scr I(E)$ is the set of all measurable  partitions  of $E$.
Obviously,  such a function $F$ is extrinsically differentiable with
  \beq\label{CL} \nn^{ext}  F(\eta)=    \sum_{i=1}^n (\pp_if)(\eta(A_1),\cdots, \eta(A_n)) \cdot 1_{A_i}.\end{equation}
We consider the square field
$$\GG_{\GA}(F,G):=  \<\GA_\eta \nn^{ext}F(\eta), \nn^{ext}G(\eta)\>_{L^2(\eta)}  = \int_E \big[\GA_\eta \nn^{ext}F(\eta)\big]\cdot \big[ \nn^{ext}G(\eta)\big]\,\d\eta,$$ and the pre-Dirichlet form
$$ \EE_{\GA,V}(F,G):= \int_{\MM} \GG_{\GA}(F,G)\, \d\BG^V,\ \ F,G\in \F C_0^\infty,$$
where $\GA$ and $V$ satisfy the following assumption.

\beg{enumerate} \item[ {\bf (H)}] For any $\eta\in \MM$, let   $\GA_\eta$ be a   bounded linear operator on $L^2(\eta)$  such that
\beq\label{H01} \<\GA_\eta h,h\>_{L^2(\eta) } \ge 0,\ \ h\in L^2(\eta),\end{equation}
  for any $A\in\B(E)$ and $x\in E$,  $\GA_\eta 1_A(x)$   is  measurable in $(\eta,x)\in \MM\times E$ and is  extrinsically   differentiable in  $\eta$  with
\beq\label{H02} \sup_{\eta(E)\le r}\big\{ \|\GA_\eta\|_{L^2(\eta)}^2 +   \|\nn^{ext}[\GA_\eta 1_A]\|_{L^2(\eta)}\big\} <\infty,\ \ r\in (0,\infty),\end{equation} where $\|\cdot\|_{L^2(\eta)} $ is the norm (or the operator norm for linear operators)  in $L^2(\eta).$

Moreover, $V\in \D(\nn^{ext})$ such that
\beq\label{H03}\sup_{\eta(E)\le r}\big\{ |V(\eta)|+ \|\nn^{ext}V(\eta)\|_{L^2(\eta)}\big\}<\infty,\ \ r\in (0,\infty).\end{equation}
\end{enumerate} Condition \eqref{H01} is essential for the nonnegativity of $\EE_{\scr A,V}$, where conditions \eqref{H02} and \eqref{H03} ensure the boundedness of  $\scr A,V$ and their extrinsic derivatives  on the level sets $\{\eta(E)\le r\}$ for $r>0$. These   conditions are standard for  establishing functional inequalities  by using perturbation argument,
see   \cite{RW04,Wbook} for the study of   the finite-dimensional models.

We write  $\GA=\1$ if $\GA_\eta$ is the identity map on $L^2(\eta)$ for every $\eta\in\MM$.
According to Theorem \ref{T2.1} below, the assumption {\bf (H)} implies that
  $(\EE_{\GA,V},\F C_0^\infty)$ is closable in $L^2(\BG^V)$ and the closure $(\EE_{\GA,V},\D(\EE_{\GA,V}))$ is a symmetric Dirichlet form. If moreover
 \beq\label{*AB2} \int_{\MM} \Big(1+\ff{\|\GA_\eta\|_{L^2(\eta)}}{1+\eta(E)}\Big)\,   \BG^V (\d\eta)<\infty,\end{equation}
 then $1\in \D(\EE_{\GA,V})$ with $\EE_{\GA,V}(1,1)=0.$  Let $(\L_{\GA,V},\D(\L_{\GA,V}))$ be the associated generator. We aim to investigate functional  inequalities for the Dirichlet form $\EE_{\GA,V}$ and the spectral gap of the generator $\L_{\GA,V}.$

We first consider  the Poincar\'e inequality
 \beq\label{PC1} \BG^V(F^2)\le \ff 1 \ll  \EE_{\GA,V}  (F,F) + \BG^V(F)^2,\ \ F\in \D(\EE_{\GA,V}),\end{equation} where $\ll>0$ is a constant.  The spectral gap of $\L_{\GA,V}$, denoted by $\gap(\L_{\GA,V})$,
  is the largest constant $\ll>0$ such that \eqref{PC1} holds.  If \eqref{PC1} is invalid, i.e. there is no any constant $\ll>0$ satisfying the inequality, we write
 $\gap(\L_{\GA,V})=0$  and   say that $\L_{\GA,V}$ does not have spectral gap. It is well known that \eqref{PC1} is equivalent to the exponential
 convergence of the associated Markov semigroup $P_t^{\GA,V}$:
 $$\|P_t^{\GA,V}F-\BG^V(F)\|_{L^2(\BG^V)}\le \e^{-\ll t}\|F\|_{L^2(\BG^V)},\ \ t\ge 0, F\in L^2(\BG^V).$$

When $\gap(\L_{\GA,V})=0$, the following weak Poincar\'e inequality was introduced in \cite{RW01}:
\beq\label{WPC}  \BG^V(F^2)\le \aa(r) \EE_{\GA,V}  (F,F) +  r\|F\|_\infty^2,\ \ F\in \D(\EE_{\GA,V}),\ \BG^V(F)=0, r>0, \end{equation}
where $\aa: (0,\infty)\to (0,\infty)$ corresponds to a non-exponential convergence rate  of $ P_t^{\GA,V} $ as $t\to\infty$, see \cite[Theorems
2.1 and 2.3]{RW01}. In particular, \eqref{WPC} implies
$$\|P_t^{\GA,V}-\BG^V\|_{L^\infty(\BG^V)\to L^2(\BG^V)}\le \inf\big\{r>0: \aa(r)\log r^{-1}\le 2t\big\}\downarrow 0\ \text{as}\ t\uparrow \infty.$$

We also consider the super Poincar\'e inequality
\beq\label{SUP} \BG^V(F^2)\le r\EE_{\GA,V}(F,F)+ \bb(r) \BG^V(|F|)^2,\ \ \ r>0, F\in \D(\EE_{\GA,V}),\end{equation}where
$\bb: (0,\infty)\to (0,\infty)$ is a decreasing function. The existence of super Poincar\'e inequality is equivalent to the uniform integrability of $P_t^{\GA,V}$ for $t>0$, and, when $P_t^{\GA,V}$ has an asymptotic density with respect to $\BG^V$, it is also equivalent to the compactness of $P_t^{\GA,V}$ in $L^2(\BG^V)$,  see \cite[Theorem 3.2.1]{Wbook} for details. According to \cite[Definition 3.1.2]{Wbook},    $P_t^{\GA,V}$  is said to have an asymptotic density, if $\|P_t^{\GA,V}-P_n\|_{L^2(\BG^V)}\to 0$ for a sequence of bounded linear operators $\{P_n\}_{n\ge 1}$ having densities with respect to $\BG^V$.
  We say that $\EE_{\GA,V}$ does not satisfy the super Poincar\'e inequality, if there is no   $\bb: (0,\infty)\to (0,\infty)$
satisfying \eqref{SUP}. In particular, \eqref{SUP} holds with $\bb(r)=\e^{cr^{-1}}$ for some constant $c>0$ if and only if the log-Sobolev inequality
  \beq\label{LSI0} \BG^V(F^2\log F^2)\le C\EE_{\GA,V}(F,F), \ \  F\in \D(\EE_{\GA,V}), \BG^V(F^2)=1\end{equation} holds for some constant $C>0$. It is well known (see \cite{BAK,GROSS}) that \eqref{LSI0} is equivalent to the hypercontractivity of $P_t^{\GA,V}$:
  $$\|P_t^{\GA,V}\|_{L^2(\BG^V)\to L^4(\BG^V)}=1\ \ \text{for\ large}\ t>0,$$ as well as the exponential convergence in entropy:
  $$\BG^V((P_t^{\GA,V} F)\log P_t^{\GA,V} F)\le \e^{- 2t/C} \BG^V(F\log F),\ \ t\ge 0,  F\ge 0, \BG^V(F)=1.$$
    See  \cite{W00a,W00b, W04} or \cite{Wbook}  for more results on the super Poincar\'e inequalities,  for instance,  estimates on the semigroup $P_t^{\GA,V}$  and higher order eigenvalues of the generator $\L_{\GA,V}$ using the function $\bb$ in \eqref{SUP}.

\

The remainder of the paper is organised as follows. In section 2, we state the main results of the paper, and illustrate these results by a typical example with specific interactions. In Section 3, we establish the integration by parts formula which implies the closability of $(\EE_{\GA,V},\F C_0^\infty)$. Then the main results are   proved in Section 4, and   extended in Section 5 to the space $\MM_s$ of finite signed measures.

\section{Main results and an example}

We first consider    $\EE_{\1,0}$ in $L^2(\BG)$ whose restriction on $\MM_1:=\{\mu\in \MM: \mu(E)=1\}$ gives rise to the Dirichlet form of the Fleming$-$Viot process.  Corresponding to results of \cite{SN,STA}  for  the Fleming$-$Viot process, we have  the   following result.
See also \cite{RW19, WZ19} for functional inequalities of different type measure-valued processes.

\beg{thm}\label{T1.2} Let $V=0$ and $\GA=\1.$ \beg{enumerate} \item[$(1)$]   $\gap (\L_{\1,0})=1,$ i.e. $\ll=1$ is the largest constant such that $\eqref{PC1}$ holds for $V=0$ and $\GA=\1$.
\item[$(2)$] If ${\rm supp}\,\nu$ contains  infinitely many points, then  $\EE_{\1,0}$ does not satisfy the super Poincar\'e   inequality.
\item[$(3)$] There exists a constant $c_0>0$ such that when ${\rm supp}\,\nu$ is a finite set,   the log-Sobolev inequality
\beq\label{LSI0'} \BG(F^2\log F^2)\le \ff{c_0}{ 1\land \dd}\,\EE_{\1,0}(F,F),\ \ F\in \D(\EE_{\1,0}), \BG(F^2)=1\end{equation}
holds, where $\dd:= \min\{\nu(\{x\}):\ x\in  {\rm supp}\,\nu\}.$  \end{enumerate}
 \end{thm}

To extend this result to $\EE_{\GA,V}$, we will adopt a split argument by making perturbations to $\EE_{\1,0}$ on bounded sets and estimating the principal  eigenvalue of $\L_{\GA,V}$ outside. To this end, we take
$$\rr(\eta)= 2\sqrt{\eta(E)},\ \ \eta\in \MM$$ and let $\BB_N=\{\eta\in\MM: \rr(\eta)\le N\}$ for $N>0$. Since  \eqref{CL} implies
\beq\label{RRD} \nn^{ext}\rr(\eta)= \ff 1 {\sqrt{\eta(E)}},\ \ \eta\in \MM\setminus\{0\},\end{equation} we have
\beq\label{*0W} \GG_{\1} (\rr,\rr):= \eta\big(|\nn^{ext} \rr(\eta)|^2\big)=\ff{\eta(E)}{\eta(E)}=1.\end{equation}

According to \eqref{LVV} below, we set
\beq\label{*TP} \L_{\GA,V} \rr(\eta)= \ff 2 {\rr(\eta)} \big[(\nu-\eta)(\GA_\eta 1) +\eta\big(\nn^{ext}[\GA_\eta1(\cdot)](\cdot)\big) +\nn^{ext}_{\GA_\eta 1}V(\eta)\big] - \ff 4 {\rr(\eta)^2} \eta(\GA_\eta 1),\end{equation} where
$$\eta\big(\nn^{ext}[\GA_\eta 1(\cdot)](\cdot)\big):=\int_{E} \nn^{ext}[\GA_\eta 1(x)](x)\,\eta(\d x).$$
Let
\beq\label{*TP0}\beg{split} &  \xi(r)= \sup_{\rr(\eta)=r} \L_{\GA,V} \rr(\eta),\ \     \underline{a}(r)= \inf_{\rr(\eta) =r}\inf_{\|\phi\|_{L^2(\eta)}=1} \<\A_\eta \phi,\phi\>_{L^2(\eta)},\\
&   \bar a(r)= \sup_{\rr(\eta) =r}\sup_{\|\phi\|_{L^2(\eta)}=1} \<\A_\eta \phi,\phi\>_{L^2(\eta)}, \ \ r>0.\end{split}\end{equation} Under {\bf(H)},  $|V(\eta)|+\|\GA_\eta\|_{L^2(\eta)}$ is bounded on $\BB_r:=\{\rr\le r\}$ for $r\in (0,\infty)$. So, these functions are bounded on $[k,K]$ for any constants $K>k>0$.  Moreover,   define
\beq\label{*TP1}    \si_k:= \sup_{t\ge k} \int_t^\infty   \e^{\int_k^r \ff{ \xi(s)}{{\underline a}(s)} \d s} \d r\int_k^t  \ff 1 {  \underline{a}(r)} \e^{-\int_k^r \ff{ \xi(s)}{\underline{a}(s)} \d s} \d r,\ \ k>0.\end{equation}
Obviously,  $ \si_k$ is non-increasing in $k$ and might be infinite. We will see in Theorem \ref{TSP}(1) that   under certain conditions $\si_k<\infty$ implies the validity of Poincar\'e inequality.

We have the following extension of Theorem \ref{T1.2} to $\EE_{\GA,V}$. When supp$\,\nu$ is finite the model reduces to finite-dimensional diffusions,   for which one may derive super Poincar\'e inequalities by making perturbations to \eqref{LSI0'}. As  the present study mainly focusses on the infinite-dimensional model,   we exclude this case in the following result.

\beg{thm}\label{TSP} Assume {\bf (H)} and   $\eqref{*AB2}$. Suppose that $  \underline{a}(r)^{-1}$ is locally bounded in $r\in [0,\infty)$ and
\beq\label{0*0} \psi(s):=   \int_0^s[\bar a (r)]^{-1/2}\,\d r\uparrow\infty \text{\  as\ } s\uparrow\infty.\end{equation}
 Then the following assertions hold. \beg{enumerate}
\item[$(1)$] If   $\lim_{k\to\infty}   \si_k<\infty\ ($equivalently, $\si_k<\infty$ for all $k>0)$, then
 $$\gap(\L_{\GA,V} )\ge \sup\bigg\{\ff 1 {2 \Phi\big(\psi^{-1}(\psi(k)+ 32  \si_k + 1 )\big) + 32   \si_k}:\ k>0\bigg\}>0,$$ where
 $$ \Phi(N):= \Big(1\lor\ff{N^2}{4\nu(E)}\Big)  \exp\Big[\sup_{\rr\le N}V -\inf_{\rr\le N}V\Big]  \sup_{r\le N}  \underline{a}(r)^{-1},\ \ N>0.$$
 \item[$(2)$]  If ${\rm supp}\,\nu$ contains  infinitely many points, then $\EE_{\GA,V}$ does not satisfy the super Poincar\'e inequality.
  \item[$(3)$] The weak Poincar\'e inequality $\eqref{WPC}$  holds for
 $$\aa(r):= \inf\Big\{2 \Phi(N):\  \BG^V(\rr>N) \le \ff r{1+r}\Big\},\ \ r>0.$$ \end{enumerate}
\end{thm}

The following result  shows that   the condition in  Theorem \ref{TSP}(1) is sharp when $\GA_\eta$ and $V(\eta)$ depend only on $\rr(\eta)$.

\beg{cor}\label{C1.2} Assume {\bf (H)} and   $\eqref{*AB2}$. Let $ V(\eta) =v(\rr(\eta))$ and $\GA_\eta =a(\rr(\eta))\1$  for  large $\rr(\eta)$ and some $a,v\in C^1([0,\infty))$ with $\underline{a}(r)>0$ for $r\ge 0$.
 Then
\beg{align*}  \xi(r):=\sup_{\rr(\eta)=r} \L_{\GA,V} \rr(\eta) = a(r)\Big(\ff{1}r+v'(r) -\ff r 2\Big) +\ff r 2   a(r),\ \ \text{for\ large}\ r> 0,\end{align*}
and   $\gap(\L_{\GA,V})>0$ if and only if $\lim_{k\to\infty}\si_k<\infty.$
 \end{cor}

As in the proof of \cite[Corollary 1.3]{RW04} using \cite[Theorem 1.1]{RW04}, it is easy to see that Theorem \ref{TSP}(2) implies the following result.

 \beg{cor}\label{C2.4} Assume {\bf (H)} and   $\eqref{*AB2}$. If $\inf_{r\ge 0}   \underline{a}(r)>0$ and  $\limsup_{ r\to\infty} \ff{ \xi(r)}{\underline{a}(r)}  <0,$
   then $\gap(\L_{\GA,V})>0.$
 \end{cor}

The above two corollaries are concerned with the validity of Poincar\'e inequality. On the other hand, according to Theorem \ref{TSP}(3), the weak Poincar\'e inequality always holds under  {\bf (H)}, \eqref{*AB2} and \eqref{0*0}. We will see in the proof that the rate function $\aa$ is derived by comparing $\EE_{\GA,V}$ with $\EE_{\1,0}$ on bounded sets $\BB_N$, $N>0$. However, when these two Dirichlet forms are far away, this $\aa$ is less sharp. As a principle, to derive a sharper weak Poincar\'e inequality, one should compare $\EE_{\GA,V}$ with a closer Dirichlet form which satisfies the Poincar\'e inequality. In this spirit, we present below an alternative result on the weak Poincar\'e inequality.
To state the result, we introduce the class $\scr H$ as follows.

\paragraph{Class   $\scr H:$} We denote $h\in \scr H$, if   $0\le h\in C^1([0,\infty))$ with $h'(r)>0$  for $r>0$, such that
\beq\label{0*1}  \xi_h(r):=  \xi(r) -\ff 2 r h(r)  \inf_{\rr(\eta)=r} \eta(\GA_\eta 1),\ \ r>0\end{equation} satisfies
\beq\label{0*2}  \si_{1,h}:= \sup_{t\ge 1} \int_t^\infty   \e^{\int_1^r \ff{ \xi_h(s)}{\underline{a}(s)} \d s} \d r\int_1^t  \ff 1 {  \underline{a}(r)} \e^{-\int_1^r \ff{ \xi_h(s)}{\underline{a}(s)} \d s} \d r<\infty.\end{equation}

\ \newline It is easy to see that $\scr H\ne \emptyset$ under the conditions of Theorem \ref{TSP} and $\inf \underline{a}>0$. For any $h\in \scr H$, let $V_h=V-h(\rr)+c(h)$, where $c(h)\in \R$ is such that $\BG^{V_h}$ is a probability measure on $\MM$.
By Theorem \ref{TSP}(1) with $k=1$, for any $h\in \scr H$,   the Poincar\'e inequality
\beq\label{PCW} \BG^{V_h}(F^2)\le C(h) \EE_{\GA,V_h}+ \BG^{V_h}(F)^2,\ \ F\in \D(\EE_{\GA,V_h})\end{equation} holds for
 \beq\label{CH} C(h):= 2\Phi_{1,h}\big(\psi^{-1}(\psi(1)+32 \si_{1,h}+1)\big)+32\si_{1,h},\ \ h\in\scr H.\end{equation}

\beg{thm}\label{NT} Assume {\bf (H)}, $\eqref{*AB2}$ and $\eqref{0*0}$. If $\scr H\ne\emptyset$,  then  $\eqref{WPC}$ holds for
  $$\aa(r):= \inf\Big\{C(h)\e^{h(N)}: \ h\in \scr H,\ N>0\ \text{with}\  \BG^V(\rr>N) \le \ff {r}{1+r}\Big\},\ \ r>0,$$
where $C(h)$ is given by $\eqref{0*2}$ and $\eqref{CH}$. \end{thm}

\

To conclude this section, we present below a simple example to illustrate the main results. For simplicity, we only consider $\GA_\eta=\1$. But  by
a simple comparison argument, the assertions apply also to $\GA_\eta$ with $\<\GA_\eta\phi,\phi\>_{L^2(\eta)}\ge c\|\phi\|_{L^2(\eta)}^2$ for some constant $c>0$ and all $\eta\in \MM$, $\phi\in L^2(\eta)$.

\beg{exa}    Consider the following potential $V_0$ with interactions given by  $\psi_i\in\B_b(E\times E), i=1,2,3:$
$$V_0(\eta)=\ff{2(\eta\times\eta)(\psi_1)}{3\eta(E)^{3/2}}+ \ff{(\eta\times\eta)(\psi_2)}{\eta(E)}+ (\eta\times\eta)(\psi_3) -p\log(1+\eta(E)),$$ where   $p\in\R$ is a constant. Let $\theta_i= \sup \psi_i, 1\le i\le 3.$ Assume that
 one of the following conditions hold:\beg{enumerate}
 \item[$(1)$] $\min\big\{\theta_3,\theta_2-1, \theta_1\cdot 1_{\{\theta_2=1\}}\big\}<0;$
 \item[$(2)$] $\theta_1=\theta_2-1=\theta_3=0$ and   $p>\nu(E)$.\end{enumerate} Then
$ Z  := \BG(\e^{V_0})\le \ff 1 {\GG(\nu(E))} \int_0^\infty  (1+s)^{-p}s^{\nu(E)-1}\e^{\theta_1s^{1/2}-(1-\theta_2)s+\theta_3s^2} \,\d s<\infty,$ so that  $\BG^V$ for  $V:=V_0-\log Z$ is a probability measure on $\MM$, and the following assertions hold:
\beg{enumerate}
\item[$(a)$]   Condition (1) implies  $\gap(\L_{\1,V})>0$;
\item[$(b)$] Under   condition $(2)$,  let
$$\theta=\max\big\{12\times 1_{\{\|\psi_3\|_\infty>0\}},\ 8\times 1_{\{\|\psi_2-1\|_\infty>0\}},\ 6\times 1_{\{\|\psi_1\|_\infty>0\}},\ 5\big\}.$$
Then there exists a constant $c>0$ such that the weak Poincar\'e inequality \eqref{WPC} holds for
       $$\aa(r)=   c r^{-\ff{\theta}{2(p-\nu(E))}},\ \ r>0. $$
 \end{enumerate}\end{exa}
\beg{proof} Obviously,  the assumptions in Theorem \ref{TSP} hold  for $V$ and $\GA_\eta=\1.$ By definition it is easy to see that
\beg{align*}  \nn^{ext}V(\eta)(x)&=  \ff{\eta(E)\eta (\psi_1(x,\cdot)+\psi_1(\cdot,x))-(\eta\times\eta)(\psi_1)}{\eta(E)^{5/2}}+\eta(\psi_3(x,\cdot)+ \psi_3(\cdot,x)) \\
 &\quad +\ff{\eta(E)\eta (\psi_2(x,\cdot)+\psi_2(\cdot,x))-(\eta\times\eta)(\psi_2)}{\eta(E)^2}
     +\ff p{1+\eta(E)}.\end{align*}
Then  \beg{align*}  \nn^{ext}_1V(\eta)  &:=\eta \big( \nn^{ext}V(\eta)\big)
 \le \theta_1\sqrt{\eta(E)} +\theta_2  \eta(E) +\theta_3\eta(E)^2+\ff{p\eta(E)}{1+\eta(E)} \\
 &= \ff{\theta_1\rr(\eta)}2   +\ff{\theta_2\rr(\eta)^2}4   + \ff{ \theta_3\rr(\eta)^4}{8}+\ff{p\rr(\eta)^2}{4+\rr(\eta)^2}.
\end{align*}

(a) If (1) holds, then $\theta_3<0,$ or $\theta_2<1$, or $\theta_3=\theta_2-1=0$ and $\theta_1<0$. In any case,   we have
$$\limsup_{\rr(\eta)\to\infty} \L_{\1,V}\rr(\eta)=\limsup_{\rr(\eta)\to\infty} \ff 2 {\rr(\eta)}\Big(\nu(E) -\ff{\rr(\eta)^2} 4 +\nn_1^{ext}V(\eta)\Big)<0,$$
so that   Corollary \ref{C2.4} implies  $\gap(\L_{\1,V})>0$.

(b) Under condition (2), we prove the weak Poincar\'e inequality for the desired $\aa(r)$. Since one may always take $\aa(r)\le 1$ in \eqref{WPC}  due to $\BG^V(F^2)\le \|F\|_\infty^2,$ it suffices to prove for small $r>0$, say $r\in (0,1].$

 It is easy to see that
\beq\label{LO1}  \BG^V(\rr>N)  \le c_0  N^{ \nu(E)-p},\ \ N>0\end{equation}   holds for some constant $c_0>0$.
 For $\vv\in (0,1]$, we take $h_\vv(s)= \vv \sqrt s$.  Since $  a=1$, it is easy to check that
$$  \si_{1,h_\vv}\le c_1 \vv^{-2}$$ for some constant $c_1>0$ independent of $\vv\in (0,1]$. Moreover, there is a constant $c_2$ independent of $\vv\in (0,1]$ such that
$$\sup_{\rr\le N} V_{h_\vv}- \inf_{\rr\le N} V_{h_\vv}\le c_2 \big[\|\psi_3\|_\infty N^4 + \|\psi_2-1\|_\infty N^2 + \|\psi_1\|_\infty N +\vv N +\log (1+N)\big].$$
Combining this with \eqref{CH}, we may find   constants $c_3,c_4>0$ independent of $\vv\in (0,1]$  such that
$$C(h_\vv)\le c_3\big(\|\psi_3\|_\infty \vv^{-12}  + \|\psi_2-1\|_\infty \vv^{-8} + \|\psi_1\|_\infty \vv^{-6} +\vv^{-5}\big)\le c_4\vv^{-\theta}.$$
Taking this into account and   applying Theorem \ref{NT} for
$$N= N_r:= \Big(\ff{2c_0   }{r}\Big)^{\ff 1 {p-\nu(E)}}, $$ such that \eqref{LO1} implies $\BG^V(\rr>N)\le \ff {r}2$ as required for $r\in (0,1]$,     we conclude that the weak Poincar\'e inequality holds for
\beg{align*}\aa(r) &:= \inf_{\vv\in (0,1]}C(h_\vv)\e^{h_\vv(N_r)}\le\inf_{\vv\in (0,1]} c_4\vv^{-\theta}
  \exp\big[\vv (2c_0r^{-1})^{\ff 1 {2(p-\nu(E))}} \big],\ \ r\in (0,1]. \end{align*} Therefore, by taking $ \vv= 1\land  r^{\ff 1 {2(p-\nu(E))}},$
we prove   \eqref{WPC} for the desired   $\aa(r).$

\end{proof}

\section{The Dirichlet form}

For any $F\in \F C_0^\infty$, let
\beq\label{LVV} \beg{split} & \L_{\GA,V} F(\eta)
 :=  \int_E \GA_\eta [\nn^{ext} F(\eta)](x)(\nu-\eta)(\d x) \\
 &\qquad  +\int_E  \nn^{ext}\big[\A_\eta(\nn^{ext}F(\eta))(x)\big](x) \, \eta(\d x)
 +  \big\<\nn^{ext} V(\eta), \GA_\eta[\nn^{ext}F(\eta)]\big\>_{L^2(\eta)}.\end{split}\end{equation}
 It is easy to see from \eqref{CL} that when $F(\eta)=f(\eta(A_1),\cdots,\eta(A_n))$ for some $n\ge 1, f\in C_0^\infty(\R^n)$ and a measurable partition
$\{A_i\}_{1\le i\le n}$ of $E$, we have
 \beg{align*} \L_{\GA,V} F(\eta)  =  &\bigg(\sum_{i=1}^n \Big[(\nu-\eta)(\GA_\eta 1_{A_i}) +\eta\big(\nn^{ext}[\GA_\eta 1_{A_i}(\cdot)](\cdot)\big) + \nn^{ext}_{\GA_\eta 1_{A_i}} V(\eta)\Big]  \pp_if \\
&  +\sum_{i,j=1}^n \eta(1_{A_i}\GA_\eta 1_{A_j})  (\pp_i\pp_j f)\bigg) (\eta(A_1),\cdots,\eta(A_n)). \end{align*}

   \beg{thm}\label{T2.1}  Assume {\bf (H)}. Then
\beq\label{*AB1} \EE_{\GA,V}(F,G) = - \int_{\MM} (G \L_{\GA,V} F) \d\BG^V,\ \ F,G\in \F C_0^\infty. \end{equation}
 Consequently, $(\EE_{\GA,V},  \F C_0^\infty)$ is closable in $L^2(\MM, \BG^V)$ whose closure $(\EE_{\GA,V},\D (\EE_{\GA,V}))$ is a symmetric Dirichlet form with generator $(\L_{\GA,V},\D(L_{\GA,V}))$ being the  Friedrichs extension of $(\L_{\GA,V}, \F C_0^\infty)$. If moreover
  $\eqref{*AB2}$ holds,    then $1\in \D(\EE_{\GA,V})$ and $\EE_{\GA,V}(1,1)=0.$
    \end{thm}

To prove this result,  we introduce the divergence operator corresponding to $\nn^{ext}$. To this end, we  formulate  the Gamma distribution  $\BG$ by using  the Poisson measure $\pi_{\hat\nu}$ with intensity $
\hat\nu(\d x,\d s):= s^{-1}\e^{-s} \nu(d x) \d s $ on $\hat E:= E\times (0,\infty)$. Recall that    $\pi_{\hat\nu}$ is the unique  probability measure on the configuration space
$$\BGG(\hat E):= \Big\{\gg=\sum_{i=1}^\infty \dd_{(x_i,s_i)}:\ \gg (K)<\infty \ \text{for\ compact\ } K\subset \hat E, (x_i,s_i)\in \hat E \Big\}$$ such that for any disjoint relatively compact subsets $\{\hat A_i\}_{1\le i\le n}$ of $\hat E$, $\{\gg\mapsto \gg(\hat A_i)\}_{1\le i\le n}$ are independent random Poisson random variables with parameters $\{\hat\nu (\hat A_i)\}_{1\le i\le n}$. Since $S(\gg):= \sum_{i=1}^\infty s_i$ for $\gg=\sum_{i=1}^\infty s_i\dd_{x_i}\in \BGG(\hat E)$ satisfies
 $$\int_{\BGG(\hat E)} S(\gg)\pi_{\hat\nu}(\d\gg) =\int_{\hat E} s\hat\nu(\d x,\d s)=\nu(E)<\infty,$$
 the measure $\pi_{\hat\nu}$ is concentrated on the   $S$-finite configuration space
 $$\BGG_{f}(\hat E):= \Big\{\gg=\sum_{i=1}^\infty \dd_{(x_i,s_i)}\in \BGG(\hat E): S(\gg):=\sum_{i=1}^\infty s_i<\infty\Big\}.$$

 \beg{lem}\label{LLK} The map
 $\Phi: \BGG_{f}(\hat E)\ni \gg=\sum_{i=1}^\infty s_i\dd_{x_i}\mapsto \sum_{i=1}^\infty s_i\dd_{x_i}\in\MM $
 is measurable with
\beq\label{*AC} \BG=\pi_{\hat\nu}\circ \Phi^{-1}.\end{equation}
Moreover,
\beq\label{MC} \beg{split} &\int_\MM \BG(\d\eta) \int_E F(\eta,x) \eta(\d x)\\
 &= \int_\MM\BG(\d\eta) \int_{\hat E} \e^{-s} F(\eta+s\dd_x, x)\nu(\d x)\d s,\ \
F\in L^1\big(\MM\times E, \BG(\d\eta) \eta(\d x)\big).\end{split}\end{equation}\end{lem}
\beg{proof} Formula \eqref{*AC} was proved in \cite[Theorem 6.2]{HKPR} for $E=\R^d$ and $\nu(\d x)=\theta\d x$ (which is an infinite measure) with $\theta>0$, by   identifying the Laplace transforms of $\BG$ and $\pi_{\hat\nu}\circ \Phi^{-1}.$ Below we explain that the same argument works to the present setting.

Firstly, the Laplace transform of $\BG$ is
\beq\label{LTR} \int_\MM \e^{-\eta(h)} \BG(\d\eta) =\e^{-\nu(\log(1+h))},\ \ h\in \B^+(E),\end{equation} where $\B^+(E)$ is the class of nonnegative measurable functions on $E$. This was given by \cite[(7)]{TVY} when $\nu$ is atomless.
In general, we decompose $\nu$ into $\nu=\nu_0+\sum_{i=1}^\infty c_i\dd_{x_i}$, where
$\nu_0$ is an atomless finite measure on $E$, $x_i\in E$ with $x_i\ne x_j$ for $i\ne j$, and  $c_i\ge 0$ with $\sum_{i=1}^\infty c_i<\infty$. Let $E_0=E\setminus\{x_i:i\ge 1\}$. By the definition of Gamma distribution,
$$\eta(h\cdot 1_{E_0}),\ \ \eta(h\cdot 1_{\{x_i\}}),\ \ \ i\ge 1$$ are independent under $\BG$,  the distribution of $\eta(h\cdot 1_{E_0})$ under $\BG$
coincides with that under $\BG_0$ (the Gamma distribution with intensity measure $\nu_0$), and the distribution of $\eta(\{x_i\})$ under $\BG$
coincides with  the one-dimensional Gamma distribution $\gg_{c_i}$ with shape parameter $c_i.$
So, applying \eqref{LTR} for $\nu_0$ replacing $\nu$ due to \cite[(7)]{TVY}, and using  the Laplace transform for Gamma distributions on $\R_+$, we derive
\beg{align*} &\int_\MM \e^{-\eta(h)} \BG(\d\eta)= \bigg(\int_\MM \e^{-\eta(h\cdot 1_{E_0})} \BG(\d\eta)\bigg)\cdot \prod_{i=1}^\infty \int_\MM \e^{-h(x_i)\eta(\{x_i\})} \BG(\d\eta)\\
&= \e^{-\nu_0(\log(1+h))}\cdot \prod_{i=1}^\infty \e^{-c_i\log (1+h(x_i))} = \e^{-\nu(\log(1+h))}.\end{align*}
Therefore, \eqref{LTR} holds.

On the other hand,   the Laplace transform for $\pi_{\hat\nu}$ (see for instance \cite{AKR}) is
$$\int_{\BGG_{pf}(\hat E)} \e^{-\gg(\hat h)}\pi_{\hat\nu}(\d\gg)= \exp\big[ -\hat\nu(1-\e^{-\hat h})\big],\ \ \hat h\in\B^+(\hat E). $$
By letting $\hat h(x,s)= s h(x)$ for $(x,s)\in\hat E$, we arrive at
\beg{align*}  & \int_{\MM} \e^{-\eta(h)}(\pi_{\hat\nu}\circ \Phi^{-1})(\d \eta)= \int_{\BGG_{pf}(\hat E) } \e^{-\gg(\hat h)} \pi_{\hat\nu} )(\d \gg)\\
&=\exp\big[-\hat\nu\big(1-\e^{-\hat h}\big)\big]= \e^{-\nu(\log(1+h))},\ \ \ h\in \B^+(E).\end{align*}
Combining this with \eqref{LTR} we prove \eqref{*AC}.

Finally, \eqref{MC} follows from \eqref{*AC} and
the Mecke formula \cite[Satz 3.1]{Mec} for Poisson measures.  \end{proof}

To establish the integration by parts formula for $\nn_\phi^{ext} F$, we introduce the divergence operator $\ddiv^{ext}$ as follows.

  Let $\phi: \MM\times E\to \R$ be measurable. If for any $x\in E$, $\phi(\cdot,x)\in \D(\nn^{ext})$  such that
$$(\BG\times\nu)(|\phi|)+ \int_{\MM} \eta\big(|\phi(\eta,\cdot)| +|\nn^{ext} \phi(\eta,\cdot)(\cdot)|\big) \BG(\d\eta)  <\infty,$$
where $\eta(\cdot)$ stands for the integral with respect to $\eta$ as in \eqref{*00}, then we write $\phi\in \D(\ddiv^{ext})$ and denote
\beq\label{DIV} \ddiv^{ext}(\phi)(\eta)= (\eta-\nu)\big(\phi(\eta,\cdot)\big) -\eta\big(\nn^{ext} \phi(\eta,\cdot)(\cdot)\big).\end{equation}

When $\phi(\eta,x)=\phi(x)$ does not depend on $\eta$, the following integration by parts formula  follows from  \cite[Theorem 14]{KLV}.   We include below a complete proof for the $\eta$-dependent $\phi$.

  \beg{lem}\label{L2.2} Let $\phi\in \D(\ddiv^{ext})$. Then
\beq\label{BBT*} \int_{\MM} (\nn^{ext}_\phi F) \, \d\BG  =  \int_{\MM} [F\ddiv^{ext}(\phi)] \,\d\BG,\ \ F\in \F C_0^\infty.\end{equation}
  \end{lem}
  \beg{proof}    By \eqref{MC} and the Dominated Convergence Theorem, we obtain
  \beg{align*} &\int_{\MM} (\nn^{ext}_\phi F) \, \d\BG= \int_{\MM\times E} \Big(\lim_{\vv\downarrow 0} \ff{F(\eta+\vv\dd_x)-F(\eta)}\vv\Big)\phi(\eta,x)\,\eta(\d x) \BG(\d\eta)\\
  &=  \int_{\MM} \BG(\d\eta) \lim_{\vv\downarrow 0}  \int_{\hat E} \ff 1 \vv \e^{-s} \big[F(\eta+(s+\vv)\dd_x)-F(\eta+s\dd_x)\big]  \phi(\eta+s\dd_x,x)\,
  \nu(\d x) \d s \\
  &=  \int_{\MM} \BG(\d\eta) \int_{\hat E}  \e^{-s}\big[\pp_s F(\eta+s\dd_x,x)\big] \phi(\eta+s\dd_x,x)\, \nu(\d x)\d s \\
  &=    \int_{\MM} \BG(\d\eta) \int_{\hat E}  \Big(\pp_s\big[\e^{-s}  F(\eta+s\dd_x)  \phi(\eta+s\dd_x,x)\big]- F(\eta+s\dd_x) \pp_s\big[
  \e^{-s} \phi(\eta+s\dd_x,x)\big]\Big) \,\nu(\d x)\d s.\end{align*}
  Noting that $F\in \F C_0^\infty$  implies $F(\eta+s\dd_x)=0$ for large $s$, we have
  $$\int_0^\infty \pp_s\big[\e^{-s}  F(\eta+s\dd_x)  \phi(\eta+s\dd_x,x)\big]\,\d s= - F(\eta) \phi(\eta,x).$$ Hence, by using \eqref{MC} again,
   \beg{align*} &\int_{\MM} (\nn^{ext}_\phi F) \, \d\BG + \int_\MM F(\eta)\nu(\phi(\eta,\cdot)) \, \BG(\d\eta)\\
  &=  -  \int_{\MM}\BG(\d\eta)    \int_{\hat E} F(\eta+s\dd_x) \e^{-s} \big[\pp_s\phi(\eta+s\dd_x,x) -\phi(\eta+s\dd_x,x)\big]\,\nu(\d x)\d s \\
  &=   \int_\MM \BG(\d\eta) \int_{\hat E} \big[\phi(\eta+s\dd_x,x)-\nn^{ext}\phi(\cdot,x)(\eta+s\dd_x)(x)\big]\e^{-s} F(\eta+ s\dd_x)\nu(\d x)\d s  \\
&= \int_\MM F(\eta) \BG(\d\eta) \int_E [\phi(\eta,x)- \nn^{ext}\phi(\eta,x)(x)] \eta(\d x).\end{align*}
Therefore, \eqref{BBT*} holds.
   \end{proof}

  \beg{proof}[Proof of Theorem \ref{T2.1}]  We first prove \eqref{*AB1}, which implies the closability of $(\EE_{\GA,V}, \F C_0^\infty)$ and that the closure is a symmetric Dirichlet form in $L^2(\BG^V)$, see \cite{FK}.
  By the definition of $\EE_{\GA,V}$ and Lemma \ref{L2.2}, for any $F,G\in\F C_0^\infty$ we have
  \beg{align*} &\EE_{\GA,V}(F,G)= \int_\MM \GG_\GA(F,G)\d\BG^V= \int_\MM \big(\nn^{ext}_{\e^V(\eta)\GA_\eta \nn^{ext} F(\eta)}G\big)(\eta)\,\BG(\d\eta)\\
  &= \int_\MM G(\eta) \ddiv^{ext}\big(\e^{V(\eta)}\GA_\eta[\nn^{ext}F(\eta)](\cdot)\big)\,\BG(\d\eta).\end{align*}
  Therefore, by \eqref{DIV}, \eqref{*AB1} holds for
  \beg{align*} \L_{\GA,V} F(\eta)&:= - \e^{-V(\eta)} \ddiv^{ext}\big(\e^{V(\eta)}\GA_\eta[\nn^{ext}F(\eta)](\cdot)\big)\\
  &=   \int_E \Big([\nn^{ext}V(\eta)(x)]\GA_\eta[\nn^{ext}F(\eta)](x) + \nn^{ext}\big(\GA_\eta[\nn^{ext}F(\eta)](x)\big)(x)\Big)\,\eta(\d x)\\
  &\quad +  \int_E  \GA_\eta[\nn^{ext}F(\eta)](x) \, (\nu-\eta)(\d x).\end{align*}

  Next, assume that  \eqref{*AB2} holds. It remains  to find   a sequence   $\{F_n\}_{n\ge 1}\subset \D(\EE_{\GA,V})$ such that
  $$\lim_{n\to\infty} \big[\BG^V(|F_n-1|^2)+\EE_{\A,V}(F_n,F_n)\big]=0.$$
  To this end, we consider
   $\rr_n:=\sqrt{n^{-1}+\rr^2},\ n\ge 1.$ By \eqref{*0W}, we have $\rr_n\in\D(\nn^{ext})$ with
 $$\GG_\1(\rr_n,\rr_n)= \ff{\rr^2}{\rr_n^2}\le 1.$$
 Let $h\in C_0^\infty([0,\infty))$ such that $h(r)=1$ for $r\le 1$ and $h(r)=0$ for $r\ge 2$. We have
 $$F_n:= h\big(n^{-1}\log[1+\rr_n ]\big)\subset \F C_0^\infty,\ \ n\ge 1.$$
 It is easy to see that  $\BG^V(|F_n-1|^2)\to 0$ as $n\to\infty$  and due to \eqref{*AB2},
 $$\limsup_{n\to\infty} \EE_{\GA,V}(F_n,F_n)\le  \limsup_{n\to\infty}  \int_{\MM}\ff{\|\GA_\eta\|_{L^2(\eta)}  \|h'\|_\infty^2}{n^2(1+\rr)^2} \,\BG^V(\d\eta)=0.$$
\end{proof}

\section{Proofs of the main results}

In this section, we prove Theorems \ref{T1.2}, \ref{TSP}, \ref{NT} and Corollary \ref{C1.2}.

 \subsection{Proof of Theorem \ref{T1.2} and a local Poincar\'e inequality}

 \beg{proof}[Proof of Theorem \ref{T1.2}] The invalidity of the super Poincar\'e inequality will be included in the proof of Theorem \ref{TSP}(3) for a more general case. So, we only prove  (1) and (3).

 (a) We first prove $\gap(\L_{\1,0})=1$, i.e. $\ll=1$ is the  optimal constant for the Poincar\'e inequality
  \beq\label{PC2} \BG(F^2)\le  \ff 1 \ll  \EE_{\1,0}  (F,F) + \BG (F)^2,\ \ F\in \F C_0^\infty \end{equation} to  hold.
 Let
 $F(\eta)= f(\eta(A_1),\cdots, \eta(A_n)) $ for some $f\in C_0^\infty(\R^n)$ and  disjoint $A_1,\cdots, A_n$.  This  Poincar\'e inequality reduces to
 $$\mu^n(f^2)-\mu^n(f)^2 \le       \mu^n \Big(\sum_{i=1}^n x_i |\pp_i f(x_1,\cdots,x_n)|^2 \Big),$$ where according to \eqref{*0},
\beq\label{*PC} \mu^n(\d x):= \prod_{i=1}^n \mu_i (\d x_i),\ \ \mu_i (\d s)=\gg_{\nu(A_i)}(\d s):=1_{[0,\infty)}(s)\ff{s^{\nu(A_i)-1}\e^{- s}}{\GG(\nu(A_i))} \d s,\ 1\le i\le n.\end{equation}  By the additive property of the Poincar\'e inequality,
 it suffices to prove that for every $1\le i\le n$, $\ll=1$ is the largest constant satisfying
$$ \mu_i(f^2)-\mu_i (f)^2 \le \ff 1 \ll \int_0^\infty r f'(r)^2 \mu_i(\d r) ,\ \  f\in C_0^\infty([0,\infty)).$$ This  follows from the fact that the generator of the Dirichlet form $$\EE_i(f,g):= \int_0^\infty r f'(r)g'(r) \mu_i(\d r),\ \ f,g\in W^{1,2}([0,\infty),\mu_i)$$  is
 $$\L_i f(r):= rf''(r)+ (\nu(A_i)-  r) f'(r),\ \ r\in [0,\infty), $$ which has spectral gap $1$ with the first eigenfunction
 $u_i(r)= r-  \nu(A_i).$

 (b) Let ${\rm supp}\,\nu=\{x_1,\cdots, x_n\}$, we have $\dd=\min\{\nu(\{x_i\}): 1\le i\le n\}>0.$ It suffices to find a universal constant $c_0>0$ such that \eqref{LSI0'} holds for
 $$F(\eta):=f(\eta(\{x_1\}),\cdots, \eta(\{x_n\})),\ \ f\in C_0^\infty(\R^n).$$
  Letting $\mu^n$ and $\mu_i$ be as in \eqref{*PC} for $A_i=\{x_i\}$,
   \eqref{LSI0'} for this $F$ becomes
$$\mu^n(f^2\log f^2)\le    \ff{c_0}{1\land\dd} \sum_{i=1}^n \int_{[0,\infty)^n}s_i(\pp_i f)^2(s_1,\cdots, s_n) \mu^n(\d s_1,\cdots, \d s_n) +\mu^n(f^2)\log\mu^n(f^2).$$
 By the additive property of the log-Sobolev inequality, this follows from the following Lemma \ref{LNN}.
  \end{proof}

  \beg{lem}\label{LNN} For any $a,b>0$, let $\mu_{a}(\d s):= 1_{[0,\infty)}(s) \ff{s^{a-1}\e^{-s}}{\GG(a)}\,\d s$ and $\mu_{a,b}(\d s):= 1_{[0,b]}(s)\ff{\mu_a(\d s)}{\mu_a([0,b])}.$ Then there exists a constant $c_0>0$ such that for any $a,b>0$,
  \beq\label{LSI'} \mu_{a,b}(f^2\log f^2) \le \ff {c_0}{a\land 1} \int_0^b s f'(s)^2 \mu_{a,b}(\d s),\ \ f\in C^1([0,b]), \mu_{a,b}(f^2)=1.\end{equation} \end{lem}

  \beg{proof} (a) Let $a\ge 2$. We will use the Bakry$-$\'Emery criterion on Riemannian manifolds with convex boundary which in particular includes $[0,b]$ for $b>0$. More precisely, let $\L_af(s)= s f''(s)+ (a-s)f'(s)$ and
  $\GG_1(f,g)(s)= sf'(s)g'(s).$ By \cite[Theorem 1.1(4)]{W10} with $\si=0$ and $t\to\infty$, if $$\GG_2 (f,f):= \ff 1 2 \L_a \GG_1(f,f)-
  \GG_1(\L_a f,f)\ge K\GG_1(f,f)$$ holds for some constant $K>0$ and all $f\in C^2([0,b])$, then
$$ \mu_{a,b}(f^2\log f^2) \le \ff {2}{K} \int_0^b s f'(s)^2 \mu_{a,b}(\d s),\ \ f\in C^1([0,b]), \mu_{a,b}(f^2)=1.$$
So, the desired inequality \eqref{LSI'} with $c_0= 4$ follows since
\beq\label{CUV} \beg{split} \GG_2(f,f) (s)&= s^2 f''(s)^2+ \ff{a+s} 2 f'(s)^2 + 2s f'(s)f''(s) \\
  &\ge \ff{a +s-2}{2 s} \GG_1(f,f)(s),\ \  s\ge 0,\end{split}\end{equation}
so that
$\GG_2(f,f)\ge \ff 1 2\GG_1(f,f)$  when $a\ge 2.$

(b) Let $a\in (0,\ff 1 2]$. By \eqref{CUV} we have $\GG_2(f,f)(s)\ge \ff{a\land 2}4 \GG_1(f,f)(s)$ for $s\ge 2.$ So, by the Bakry$-$\'Emery criterion,
\beq\label{LSI2}\mu_{a,b_1}(1_{[2,b_1]} f^2\log  f^2)\le \ff 8{a\land 2} \mu_{a,b_1}(1_{[2,b_1]}\GG_1(f,f))+ \mu_{a,b_1}(1_{[2,b_1]} f^2)   \log \mu_{a,b_1}(1_{[2,b_1]} f^2)\end{equation} holds for any $b_1>2$ and all $f\in C^1([0,b_1]).$

On the other hand, for any $b_2>0$ and $f\in C^1([0,b_2])$ with $\mu_{a,b_2}(f)=0$, there exists $r_0\in [0,b_2]$ such that $f(r_0)=0.$ So, for any $r\in [0, b_2]$ we have
\beg{align*} |f(r)|&=\bigg|\int_{r_0}^r f'(s)\d s\bigg|\le \bigg(\int_0^{b_2} s f'(s)^2\mu_{a,b_2}(\d s)\bigg)^{\ff 1 2}\bigg(\int_0^{b_2} s^{-a}\e^{s} \GG(a)\,\d s\bigg)^{\ff 1 2}\\
&\le \Big(\ff{\GG(a)b_2^{1-a}\e^{b_2}}{1-a}  \mu_{a,b_2}(\GG_1(f,f)) \Big)^{\ff 1 2},\ \ r\in [0, b_2].\end{align*}
Therefore, for  $\mu_{a,b_2}(f^2)=1$ with $\mu_{a,b_2}(f)=0$ we have
\beg{align*} \mu_{a,b_2}(f^2\log f^2)&\le \mu_{a,b_2}(f^2) \log\Big[\ff{\GG(a)b_2^{1-a}\e^{b_2}}{1-a} \mu_{a,b_2}(\GG_1(f,f))\Big]\\
&\le  \ff{\GG(a)b_2^{1-a}\e^{b_2}}{1-a} \mu_{a,b_2}(\GG_1(f,f))-1.\end{align*}This implies
\beq\label{AAA} \beg{split} &\mu_{a,b_2}(f^2\log f^2) -\mu_{a,b_2}(f^2)\log \mu_{a,b_2}(f^2)\\
&\le  \ff{\GG(a)b_2^{1-a}\e^{b_2}}{1-a} \mu_{a,b_2}(\GG_1(f,f))-\mu_{a,b_2}(f^2),\ \ f\in C^1([0,b_2]),\ \mu_{a,b_2}(f)=0.\end{split} \end{equation}
In general, for a non-zero function   $f\in C^1([0,b_2])$, let $\tt f= f -\mu_{a,b_2}(f)$. We have (see \cite{BAK})
\beq\label{CCC} \beg{split} &\mu_{a,b_2}(f^2\log f^2)- \mu_{a,b_2}(f^2)\log \mu_{a,b_2}(f^2)\\
 &\le \mu_{a,b_2}(\tt f^2\log \tt f^2)- \mu_{a,b_2}(\tt f^2)\log \mu_{a,b_2}(\tt f^2) + 2 \mu_{a,b_2}(\tt f^2).\end{split} \end{equation} Combining this with \eqref{AAA} and using the Poincar\'e inequality \eqref{TT1} below, we arrive at
\beq\label{BBB} \beg{split} & \mu_{a,b_2}(f^2\log f^2)- \mu_{a,b_2}(f^2)\log \mu_{a,b_2}(f^2)\\
&\le \Big(\ff{\GG(a)b_2^{1-a}\e^{b_2}}{1-a} +1\Big) \mu_{a,b_2}(\GG_1(f,f)),\ \ b_2>0, f\in C_b^1([0,b_2]).\end{split}\end{equation}

In conclusion, when $b\le 4$, the desired inequality \eqref{LSI'} for $a\in (0,\ff 1 2]$ follows from \eqref{BBB}. Finally,    for $b\ge 4$ we deduce from \eqref{LSI2} and \eqref{BBB} that for any $f\in C^1([0,b])$ with $\mu_{a,b}(f^2)=1$,
\beg{align*} &\mu_{a,b}(f^2\log f^2) =\ff{\int_0^2s^{a-1}\e^{-s}\d s}{\int_0^b s^{a-1}\e^{-s}\d s} \mu_{a,2}(f^2\log f^2) +  \mu_{a,b}(1_{[2,b]}f^2\log f^2)\\
&\le \Big(\ff{\GG(a)2^{1-a}\e^{2}}{1-a} +1\Big) \mu_{a,b}(1_{[0,2]}\GG_1(f,f))+ \mu_{a,b}(1_{[0,2]} f^2)   \log \ff {\GG(a)} {\int_0^2 s^{a-1} \e^{-s}\d s}\\
&\quad +\ff 8{a\land 2} \mu_{a,b}(1_{[2,b]}\GG_1(f,f))+ \mu_{a,b}(1_{[2,b]} f^2)   \log \ff {\GG(a)} {\int_2^b s^{a-1} \e^{-s}\d s}\\
&   \le \ff{c_1}{a} \int_0^b s f'(s)^2\,\mu_{a,b}(\d s) + \ff{c_1}a \mu_{a,b}(f^2),\end{align*} where
$c_1>0$ is a universal constant independent of $a\in (0,\ff 1 2]$ and $b\ge 4.$ Combining this with \eqref{CCC} and the Poincar\'e inequality \eqref{TT1} below, we prove the inequality \eqref{LSI'} for some universal constant $c_0>0$ and all $a\in (0,\ff 1 2]$ and $b\ge 4. $

(c) Let $a\in (\ff 1 2, 2)$. In this case, we have $a':=\ff a 4\in (0,\ff 1 2]$, so that by (b)  there exists a constant $c_0>0$ such that
\beq\label{GGG} \mu_{a',b}(f^2\log f^2) \le \ff {c_0}{a}\int_0^b s f'(s)^2 \mu_{a',b}(\d s),\ \ \ a\in \Big(\ff 1 2,2\Big), f\in C^1([0,b]), \mu_{a',b}(f^2)=1.\end{equation} Let $\bar \mu_{a',\infty}(\d s_1,\d s_2,\d s_3, \d s_4)=\prod_{i=1}^4 \mu_{a',\infty}(\d s_i),$ where $\mu_{a',\infty}:=\lim_{b\to\infty} \mu_{a',b}$ is the Gamma distribution with parameter $a'$. By the property of Gamma distributions we have
$$\int_{[0,\infty)^n} f(s_1+s_2+s_3+s_4) \bar\mu_{a',\infty}(\d s_1,\d s_2,\d s_3, \d s_4)=\int_{[0,\infty)}f(s)\mu_{a,\infty}(\d s),\ \ f\in \B_b([0,\infty)).$$
Using \eqref{GGG} with $b\to\infty$ and the additivity property of the log-Sobolev inequality, we obtain
\beg{align*} &\bar\mu_{a',\infty}( F^2\log F^2)-  \bar \mu_{a',\infty}( F^2) \log \bar \mu_{a',\infty}( F^2) \\
&\le  \ff {c_0}{a}\int_0^b \sum_{i=1}^4 s_i \pp_i F (s_1,\cdots , s_4)^2 \bar \mu_{a',b}(\d s_1,\cdots, \d s_4),
\ \ F\in C^1_b([0,\infty)^4).\end{align*}
By an approximation argument we may apply this inequality to  $$F(s_1,\cdots, s_4):= f(b\land (s_1+\cdots+s_4))$$ for $f\in C^1([0,b]),$ so that  \eqref{LSI'} is derived.
\end{proof}

To prove Theorem \ref{TSP}, we   consider the local Poincar\'e   inequality for $\EE_{\1,0}$ on the set $\BB_N$,  by  decomposing $\eta$ into the radial part $\eta(E)$ and the simplicial  part $ \bar\eta:= \ff{\eta}{\eta(E)}$. It is well known that under $\BG$ these two parts are independent with
\beq\label{OPQ} \beg{split} & \BG(\eta(E)<r, \bar \eta\in \mathbf A)= {\bf Dir}(\mathbf A) \gg_{\nu(E)}([0,r)),\ \ r>0,  \mathbf A\in \B(\MM_1),\end{split}\end{equation}
where $\gg_{\nu(E)}(\d s):=  1_{[0,\infty)}(s) \ff{s^{\nu(E)-1}\e^{-s}}{\GG(\nu(E))}\,\d s, $ and {\bf Dir}  is the Dirichlet measure with intensity measure $\nu$, see for instance \cite{STA} for details. According to \cite{SN} (see also \cite[Proposition 3.3]{STA}), we have the Poincar\'e inequality
\beq\label{PDI} {\bf Dir}(F^2)\le {\bf Dir}(\GG^D(F,F))+ {\bf Dir}(F)^2,\ \ F\in \F C_0^\infty,\end{equation}
where for $F(\eta)= f(\eta(A_1),\cdots, \eta(A_n))$ and $\eta\in\MM_1$,
\beq\label{GGD} \GG^D(F,F)(\eta):= \sum_{i,j=1}^n \big[\dd_{ij} \eta(A_i)-\eta(A_i)\eta(A_j)\big]\cdot[(\pp_i f)(\pp_j f)](\eta(A_1),\cdots,\eta(A_n)).\end{equation}

 \beg{lem} \label{CN}  For any $N>0$,
 \beq\label{DIR2} \BG(1_{\BB_N}F^2)\le \Big(\ff{N^2}{4\nu(E)}\lor 1\Big)\BG(1_{\BB_N}\GG_{\1}(F,F)),\ \ F\in \F C_0^\infty, \BG(1_{\BB_N}F)=0.\end{equation}
 \end{lem}

\beg{proof}  Since $\BB_N=\{\eta(E)\le N^2/4\}$,   \eqref{OPQ} implies
\beq\label{TT0} \int_{\MM} \big[1_{\BB_N} F\big](\eta)\BG(\d\eta) = \int_{\MM_1\times [0,N^2/4]} F(s\bar \eta) {\bf Dir}(\d\bar \eta) \gg_{\nu(E)}(\d s),\ \ F\in L^1(1_{\BB_N}\BG).\end{equation}
We   observe that \eqref{T1.2} implies
\beq\label{TT1} \gg_{\nu(E)}(1_{[0,r]}f^2)\le\int_0^r sf'(s)^2 \gg_{\nu(E)}(\d s),\ \ r>0, f\in C^1([0,r]), \gg_{\nu(E)}(1_{[0,r]}f)=0.\end{equation}
Indeed, applying the Poincar\'e inequality
$$\BG(F^2)\le \EE_{\1,0}(F,F) + \BG(F)^2$$ to
$F(\eta):= f(\eta(E)\land r),$ and noting that for $\tt f(s):=f(s\land r)$ we have
\beg{align*} &\BG(F^i)= \gg_{\nu(E)}(\tt f^i)= \gg_{\nu(E)}(1_{[0,r]}f^i)+ \gg_{\nu(E)}((r,\infty))f(r),\ i=1,2,\\
& \EE_{\1,0}(F,F) = \int_0^\infty s \tt f'(s)^2\d s= \int_0^N s   f'(s)^2\d s,\end{align*}
 it follows that
\beg{align*}&\gg_{\nu(E)}(1_{[0,r]}f^2)= \gg(\tt f^2) - \gg_{\nu(E)}((r,\infty))f(r)^2 \\
&\le \int_0^N s   f'(s)^2\d s + \gg_{\nu(E)}((r,\infty))^2 f(r)^2 -\gg_{\nu(E)}((r,\infty))f(r)^2 \le\int_0^N s   f'(s)^2\d s.\end{align*}
By the additivity property of the Poincar\'e inequality,  \eqref{PDI}, \eqref{TT0} and \eqref{TT1}, we obtain that for any
 $F\in \F C_0^\infty$ with $ \BG(1_{\BB_N} F)=0$,
\beg{align*}  \BG(1_{\BB_N} F^2) &\le  \int_{\MM_1\times [0,N^2/4]}\Big[ \ff 1 {\nu(E)} \GG^D(F(s\cdot), F(s\cdot))(\bar \eta) + s \Big|\ff{\pp}{\pp s} F(s\bar \eta)\Big|^2 \Big]{\bf Dir}(\d\bar \eta)\gg_{\nu(E)}(\d s)\\
& = \int_{\BB_N} \Big[ \ff 1 {\nu(E)} \GG^D(F(\eta(E)\cdot), F(\eta(E)\cdot))(\bar \eta) + \eta(E) \Big|\ff{\pp}{\pp \eta(E)} F(\eta(E)\bar \eta)\Big|^2 \Big]\, \BG(\d\eta).\end{align*}
So, it remains to prove
\beq\label{TPP} \beg{split} I(\eta)&:=\ff 1 {\nu(E)} \GG^D(F(\eta(E)\cdot), F(\eta(E)\cdot))(\bar \eta) + \eta(E) \Big|\ff{\pp}{\pp \eta(E)} F(\eta(E)\bar \eta)\Big|^2\\
& \le \Big(\ff{N^2}{4\nu(E)} \lor 1\Big) \GG_{\1}(F,F)(\eta),\ \ \eta(E)\le \ff {N^2}4.\end{split} \end{equation}
For $F\in \F C_0^\infty$ with  $F(\eta)=f\big(\eta(A_1),\cdots, f(A_n))= f(\eta(E)\bar \eta(A_1),\cdots, \eta(E) \bar\eta(A_n)\big),$
by \eqref{GGD} we have
\beg{align*} &\GG^D\big(F(\eta(E)\cdot),F(\eta(E)\cdot)\big)(\bar \eta)\\
&= \sum_{i,j=1}^n \big[\dd_{ij}  \bar\eta(A_i) - \bar\eta(A_i) \bar\eta(A_j)\big] \eta(E)^2 [(\pp_i f)(\pp_j f)](\eta(A_1),\cdots,\eta(A_n))\\
&= \sum_{i,j=1}^n \big[\dd_{ij} \eta(A_i)\eta(E)-\eta(A_i)\eta(A_j)) [(\pp_i f)(\pp_j f)](\eta(A_1),\cdots,\eta(A_n)).\end{align*}
Moreover,
\beg{align*} \eta(E)\Big|\ff{\pp}{\pp\eta(E)} F(\eta(E) \bar\eta)\Big|^2&= \eta(E)\bigg|\sum_{i=1}^n  \bar\eta(A_i) (\pp_i f) (\eta(A_1),\cdots, \eta(A_n))\bigg|^2 \\
&= \ff 1 {\eta(E)} \sum_{i,j=1}^n \eta(A_i)\eta(A_j) [(\pp_i f)(\pp_j f)](\eta(A_1),\cdots,\eta(A_n)).\end{align*}
So, when $\eta(E)\le \ff {N^2} 4$ (i.e. $\rr(E)\le N$),
\beg{align*}I(\eta)&\le \ff{\eta(E)\lor \nu(E)}{\nu(E)} \Big(\ff 1 {\eta(E)} \GG^D(F(\eta(E)\cdot), F(\eta(E)\cdot))( \bar\eta) + \eta(E) \Big|\ff{\pp}{\pp \eta(E)} F(\eta(E) \bar\eta)\Big|^2\Big)\\
&= \Big(1\lor\ff{N^2}{4\nu(E)}\Big)\GG_{\1}(F,F)(\eta).\end{align*}
This implies \eqref{TPP},  and hence finishes the proof.
\end{proof}

 \subsection{Proofs of   Theorem \ref{TSP} and Corollary \ref{C1.2}}

   \beg{proof}[Proof of Theorem \ref{TSP}]  We will make a standard split argument by using the local Poincar\'e  inequality  \eqref{DIR2} and the principal  eigenvalue of $\L_{\GA,V}$  outside $\BB_N$.  To estimate the principal eigenvalue,
we     recall   Hardy's criterion for the first mixed eigenvalue. Consider the following differential operator on $[0,\infty)$:
$$  \L f(r) =   \underline{a}(r)f''(r)+  \gg(r)f'(r),\ \ r\ge 0.$$ For any $k>0$ and $n\ge 1$, let   $ \ll_{k,n}$ be the first mixed eigenvalue of $ \L$ on $[k,k+n]$ with Dirichlet boundary condition at $k$ and Neumann boundary condition at $k+n$. Define
$$ \si_{k,n}= \sup_{t\in(k,n+ k)} \int_t^{n+k}\e^{\int_k^r \ff{ \gg(s)}{\underline{a}(s)} \d s}  \d r\int_k^t  \ff 1 {  \underline{a}(r)} \e^{-\int_k^r\ff{ \gg(s)}{\underline{a}(s)} \d s}  \d r.$$ By Hardy's criterion, see for instance \cite[Theorem 1.4.2]{Wbook}, we have
\beq\label{*A1} \ff 1 { \si_{k,n}}\ge  \ll_{k,k+n} \ge  \ff 1 {4 \si_{k,n}},\ \ n\ge 1, k>0.\end{equation}
Below we prove assertions (1)-(3) respectively.

  (1) By \eqref{DIR2} and a standard perturbation argument, we have
 \beq\label{*A5} \BG^V(1_{\{\rr\le N\}}F^2)\le \BG^V(1_{\{\rr\le N\}}F)^2 +    \Phi(N) \BG^V(1_{\{\rr\le N\}}\GG_\GA(F,F)),\ \ F\in \F C_0^2.\end{equation}
 If $  \si_k<\infty$ for some $k>0$,  it suffices  to prove the Poincar\'e inequality
 \beq\label{*A6}  \BG^V(F^2)\le   C  \BG^V( \GG_\GA(F,F)),\ \ F\in \F C_0^2, \BG^V(F)=0\end{equation}  for
 $$C= 2 \Phi\big(\psi^{-1}\Big(\psi(k)+ 8  \ll_k^{-1}+1\big)+ \ll_k\Big) +    8   \ll_k^{-1},$$where according to \eqref{*A1},
\beq\label{PPK0}  \ll_k:=\lim_{n\to\infty} \ll_{k,n} \ge \ff 1{4 \si_k}.\end{equation}
 Let $F\in \F C_0^2$ such that ${\rm supp}\, F\subset\BB_{N_1} $ for some constant $N_1>k$. For any $N\ge k$, let
 $$F_N= F[(\psi(\rr)-\psi(N))^+\land 1].$$
 Then $F_N=0$ for $\rr\le N$ and $F_N=F$ for $\psi(\rr)\ge \psi(N)+1$. For $n>N_1$, let $u_{n}\ge 0$ be the first mixed eigenfunction of $ \L$ on $[k, k+n]$ with Dirichlet boundary condition at $k$ and Neumann boundary condition at $k+n$, such that
 $$ u_n(k)= u_n'(k+n)=0,\ \  u_n'(r)> 0\ \text{for}\ r\in (k,k+n),\ \  \L u_n=- \ll_{k,n}u_n\le 0.$$
Combining this with the definition of $ \L$ we obtain
 $$\L_{\GA,V} (u_n\circ\rr)\ge ( \L u_n)\circ\rr,\ \ \rr\in [k,k+n].$$ So,
  \beq\label{*CD1} \beg{split}  & \ll_{k,n} \BG^V(F_N^2)
 = - \int_{\{k<\rr<n+k\}} \ff{F_N^2}{u_n^2\circ\rr} (- \L u_n)\circ \rr\, \d\BG^V \\
 &\le  - \int_{\{k<\rr<n+k\}} \ff{F_N^2}{u_n\circ\rr} \big[-\L_{\GA,V} ( u_n\circ \rr)\big]  \d\BG^V.\end{split}\end{equation}
 To apply the integration by parts formula, we   approximate  $u_n$ as follows.  Since  $u_n(k)=u_n'(k+n)=0$, we may construct a sequence $\{u_{n,m}\}_{m\ge 1} \subset C^\infty([0,\infty))$ such that
\beg{align*} &u_{n,m}(r)=u_n(r)\ \text{for}\ r\in [k+m^{-1}, k+n-m^{-1}],\\
& u_{n,m}(r)=0\ \text{for}\ r\le k,\ \ u_{n,m}'(r)=0\ \text{for }\ r\ge k+n,\\
& \sup_{m\ge 1} \sup_{r\ge k} \big(|u_{n,m}'(r)|+|u_{n,m}''(r)|\big)<\infty.\end{align*}
Since $F_N=0$ for $\rr\le N$, \eqref{*CD1} implies that for any $k<N$,
\beq\label{*CD2} \beg{split}   \ll_{k,n} \BG^V(F_N^2) &=-\lim_{m\to\infty} \int_{\{k<\rr<n+k\}} \ff{F_N^2}{u_{n,m}^2\circ\rr} (- \L u_{n,m})\circ \rr\, \d\BG^V\\
&= \lim_{m\to\infty} \int_\MM \Big\<\GA_\eta \nn^{ext} \ff{F_N^2}{u_{n,m}\circ \rr}(\eta),  \nn^{ext} (u_{n,m}\circ\rr)(\eta)\Big\>_{L^2(\eta)} \d\BG^V.\end{split}\end{equation}
On the other hand, since $\GA_\eta$ is positive definite due to {\bf (H)}, for any $u\in C^2([0,\infty))$ with $u(r)>0$ for $r\ge N$, we have
\beg{align*} & \Big\<\GA_\eta \nn^{ext} \ff{F_N^2}{u\circ \rr}(\eta),  \nn^{ext} (u\circ\rr)(\eta)\Big\>_{L^2(\eta)} \\
 &=  \big\<\GA_\eta\nn^{ext} F_N(\eta), \nn^{ext} F_N(\eta)\big\>_{L^2(\eta)} \\
     &\qquad - \Big\<\GA_\eta\Big[\nn^{ext} F_N -\ff{F_N}{u\circ \rr} \nn^{ext}(u\circ\rr)\Big](\eta), \nn^{ext} F_N (\eta)-\ff{F_N}{u\circ \rr} \nn^{ext}(u\circ\rr)(\eta)\Big>_{L^2(\eta)}   \\
 &\le    \big\<\GA_\eta\nn^{ext} F_N(\eta), \nn^{ext} F_N(\eta)\big\>_{L^2(\eta)}.\end{align*}
Combining this with  \eqref{*CD2} and the definition of $F_N$, we obtain
\beg{align*} & \ll_{k,n} \BG^V(F_N^2)\le  \int_\MM \<\GA_\eta\nn^{ext}F_N,\nn^{ext} F_N\>_{L^2(\eta)}   \d\BG^V\\
&\le 2 \EE_{\GA,V}(F,F) + 2  \int_{\{\psi(N)<\psi(\rr)< \psi(N)+1\}}F^2 \GG_\GA(\psi(\rr),\psi(\rr)) \d\BG^V.\end{align*} Multiplying  by $ \ll_{k,n}^{-1} $ and letting $n\to\infty$ leads to
 \beq\label{IMP} \int_{\MM} F_N^2 \d \BG^V\le \ff 2 { \ll_k} \EE_{\GA,V}(F,F) + \ff 2 { \ll_k} \int_{\{\psi(N)<\psi(\rr)< \psi(N)+1\}}F^2 \GG_\GA(\psi(\rr),\psi(\rr)) \d\BG^V.\end{equation}
 By the definition of $\psi$ and $\bar a$, and noting that $\GG_\1(\rr,\rr)=1$, we have
 \beq\label{GMM} \GG_\GA(\psi(\rr),\psi(\rr))(\eta)= \ff{\<\GA_\eta\nn^{ext}\rr(\eta), \nn^{ext}\rr(\eta)\>_{L^2(\eta)}}{\bar a (\rr(\eta))}
 \le \GG_\1(\rr,\rr)(\eta)=1.\end{equation}
So, \eqref{IMP} implies
\beq\label{IMP'}   \int_{\MM} F_N^2 \d \BG^V\le   \ff 2 { \ll_k} \EE_{\GA,V}(F,F) + \ff 2 { \ll_k} \int_{\{\psi(N)<\psi(\rr)< \psi(N)+1\}}F^2 \d\BG^V.\end{equation}
Letting $\[s\]=\sup\{k\in\mathbb Z:\ k\le s\}$ be the integer part of a real number $s$, we have
 $$\int_\MM F^2\d\BG^V\ge \sum_{i=1}^{1+ \[8 \ll_k^{-1}\]}\int_{\{\psi(k)+ i-1  <\psi(\rr)< \psi(k) +i ]} F^2\d\BG^V.$$
 Then there exists $N\in \big[k, \psi^{-1}\big(\psi(k)+ 8  \ll_k^{-1}\big)\big]$ such that
 $$\int_{\{\psi(N)<\psi(\rr)< \psi(N)+1\}}F^2 \d\BG^V\le \ff{ \ll_k} 8 \int_{\MM} F^2 \d\BG^V,$$ so that \eqref{IMP'} yields
 \beq\label{PPK} \int_{\MM} F_N^2 \d \BG^V\le \ff 2 { \ll_k} \EE_{\GA,V}(F,F) + \ff 1 4 \int_{\MM} F^2 \d\BG^V.\end{equation}
 Combining this with \eqref{*A5} and noting that $\BG^V(F)=0$, we may find    $N\in [k,\psi^{-1}(\psi(k)+  8   \ll_k^{-1})]$ such that
 \beg{align*} &\int_{\MM} F^2 \d\BG^V\le \int_{\psi(\rr)\le \psi(N)+1}  F^2 \d\BG^V + \int_{\MM} F_N^2 \d\BG^V\\
 &\le  \Phi(\psi^{-1}(\psi(N)+1)) \EE_{\GA,V}(F,F) + \BG^V(1_{\{\psi(\rr)\ge \psi(N)+1\}}F)^2 +  \int_{\MM} F_N^2 \d\BG^V \\
 &\le      \Phi(\psi^{-1}(\psi(N)+1)) \EE_{\GA,V}(F,F) + 2\int_{\MM} F_N^2 \d\BG^V\\
 &\le \Big(    \Phi(\psi^{-1}(\psi(N)+1)) +\ff 4 { \ll_k}\Big) \EE_{\GA,V}(F,F)+ \ff 1 2 \int_{\MM} F^2 \d\BG^V.\end{align*}
Since $\Phi(N)$ is increasing in $N\in \big[k, \psi^{-1}\big(\psi(k)+ 8   \ll_k^{-1}\big)\big]$,  this implies \eqref{*A6}  with
$$C=2\Phi(\psi^{-1}(\psi(N)+1) )+ \ff 8 { \ll_k} \le 2  \Phi\big(\psi^{-1}\Big(\psi(k)+ 8  \ll_k^{-1}+1\big)+ \ll_k\Big) +    8   \ll_k^{-1}. $$   Then the proof is finished by \eqref{PPK0}.

 (2) Assume that supp\,$\nu$ is an infinite set.  To disprove the super Poincar\'e inequality, it suffices to construct a sequence $\{F_n\}\subset \D(\EE_{\GA,V})$ such that
$\BG^V(F_n^2)>0$ and
\beq\label{AAO} C:= \sup_{n\ge 1} \ff{\EE_{\GA,V}(F_n,F_n)}{\BG^V(F_n^2)}<\infty,\ \ \lim_{n\to\infty} \ff{\BG^V(|F_n|)^2}{\BG^V(F_n^2)}=0.\end{equation}
Indeed, if \eqref{SUP} holds for some $\bb: (0,\infty)\to (0,\infty)$, then
$$1\le r \ff{\EE_{\GA,V}(F_,F_n)}{\BG^V(F_n^2)}+\bb(r) \ff{\BG^V(|F_n|)^2}{\BG^V(F_n^2)},\ \ n\ge 1, r>0.$$
Combining this with \eqref{AAO} and letting $n\to\infty$, we obtain
$1\le rC $ for all $r>0$ which is impossible.

We now show that \eqref{AAO} holds for $F_n(\eta):= (1-\eta(E))^+ \ff{\eta(A_n)}{\eta(E)}$, where
  $\{A_n\}_{n\ge 1}$ are measurable subsets of $E$ such that
$\ff 1 2 \nu(E)> p_n:= \nu(A_n)\downarrow 0$ as $n\uparrow\infty$, which exist since supp\ $\,\nu$ is an infinite set.

Obviously, $\{F_n\}_{n\ge 1}\subset \D(\EE_{\GA,V}).$
Since $\|\GA_\eta\|_{L^2(\eta)} + \e^{V(\eta)}+\e^{-V(\eta)}$ is bounded  on the set $\{\eta(E)\le 1\}$, we may find   constants $K_i, C_i>0, i=1,2,3$ such that   \eqref{OPQ} implies for all $n\ge 1$ that
 \beg{align*} &\BG^V(F_n^2)\ge K_1 \BG(F_n^2) = K_1 \int_0^1 \ff{(1-s)^2s^{\nu(E)-1}\e^{-s}}{\GG(\nu(E))}\,\d s\int_0^1 \ff{t^{p_n+1} (1-t)^{\nu(E)-p_n-1}}{\GG(p_n)\GG(\nu(E)-p_n)}\,\d t\ge C_1 p_n,\\
 &\BG^V(|F_n|)^2\le K_2 \BG(|F_n|)^2 = K_2 \bigg(\int_0^1 \ff{(1-s)s^{\nu(E)-1}\e^{-s}}{\GG(\nu(E))}\,\d s\int_0^1 \ff{t^{p_n} (1-t)^{\nu(E)-p_n-1}}{\GG(p_n)\GG(\nu(E)-p_n)}\,\d t\bigg)^2\le C_2 p_n^2,\\
  &\EE_{\GA,V}(F_n,F_n)\le K_3\BG(\|\nn^{ext}F_n\|^2_{L^2(\eta)} )= K_3 \int_{\{\eta(E)\le 1\}} \BG(\d\eta) \int_E |(1-\eta(E))- \eta(A_n)|^2 \d\eta\\
   &\qquad \le 2 K_3 \int_{\{\eta(E)\le 1\}} \big[(1-\eta(E))^2 \eta(A_n) +\eta(A_n)^2\eta(E)\big]\BG(\d\eta)\\
   &\qquad \le 4 K_3 \int_{\{\eta(E)\le 1\}} \ff{\eta(A_n)}{\eta(E)}\, \BG(\d\eta)\le C_3 p_n. \end{align*} Since $p_n\downarrow 0$ as $n\uparrow \infty$, we prove \eqref{AAO}.

  (3) The local Poincar\'e inequality \eqref{DIR2} implies that for any $F\in \F C_0^\infty$ with $\BG^V(F)=0$,
 \beg{align*}\BG^V(F^2) &= \BG^V(1_{\BB_N} F^2)+ \BG^V(F^2\cdot 1_{\BB_N^c})\\
 &\le  2   \Phi(N)\EE_{\GA,V}(F,F) + \ff{1}{\BG^V(\BB_N)} \BG^V(1_{\BB_N^c} F)^2 + \BG^V(\rr>N) \|F\|_\infty^2\\
 &\le  2  \Phi(N) \EE_{\GA,V}(F,F)+ \Big(\ff{\BG^V(\rr>N)^2}{1-\BG^V(\rr>N)}+\BG^V(\rr>N)\Big) \|F\|_\infty^2,\ \ N>0.\end{align*}
 So, for any $r>0$, taking $N>0$ such that $\ff{\BG^V(\rr>N)^2}{1-\BG^V(\rr>N)}+\BG^V(\rr>N)\le r$, i.e. $\BG^V(\rr>N)\le \ff r {1+r}$, we prove \eqref{WPC}.

\end{proof}

\beg{proof}[Proof of Corollary \ref{C1.2}]
Let $r_0\in (0,\infty)$ such that $\GA_\eta =a(\rr(\eta))\1$ and $V(\eta)=v(\rr(\eta))$ for large $\rr(\eta)\ge r_0.$ By \eqref{RRD}   we have
 \beg{align*} & (\nu-\eta)(\GA_\eta 1)= a(\rr(\eta)) \big(\nu(E)- \eta(E)\big)= a(\rr(\eta)) \Big(\nu(E)-\ff{\rr(\eta)^2}4\Big),\\
 &\eta\big(\nn^{ext}[\GA_\eta 1 (\cdot)](\cdot)\big)= a'(\rr(\eta)) \eta\Big(\ff 1 {\sqrt{\eta(E)}}\Big) =\ff {\rr(\eta)} 2 a'(\rr(\eta)),\\
 &\nn^{ext}_{\GA_\eta 1} V(\eta):=\eta\big([\GA_\eta 1] \nn^{ext} V(\eta)\big)= (av')(\rr(\eta)) \ff {2\eta(E)} {\rr(\eta)} = \ff {\rr(\eta)} 2 (av')(\rr(\eta)),\\
 & \ff 4 {\rr(\eta)^3} \eta(\GA_\eta 1)= \ff{4 a(\rr(\eta))\eta(E)}{\rr(\eta)^3}= \ff{a(\rr(\eta))}{\rr(\eta)},\ \ \rr(\eta)\ge r_0. \end{align*}
 Then \eqref{*TP} implies
\beq\label{PKKN} \L_{\GA_\eta,V} \rr(\eta)= \xi(\rr(\eta))\end{equation}  for the given function $\xi$. So, when $\si_k<\infty$ for some $k>0$, Theorem \ref{TSP}(1) implies
$\gap(\L_{\GA,V})>0.$

On the other hand, let $ \si_k=\infty$ for all $k>0$. We have
\beq\label{**}  \ll_k:=\lim_{n\to\infty}\ll_{k,n}=0,\ \ k>0,\end{equation} where $\ll_{k,n}$ is given in the proof of Theorem \ref{TSP}.
Let $u_{k,n}$ be the corresponding   first mixed eigenfunction of $\L$ on $[k,k+n]$ with $u_{k,n}(r)>0$ in $(k,k+n]$, and let
$$\Theta_v(\d s) = \ff{\e^{v(s)-s}s^{\nu(E)-1}}{\GG(\nu(E))}\,\d s,$$ such that $\L$ is symmetric in $L^2([k,k+n],\Theta_v)$ under the mixed boundary conditions. Then
$$\int_k^{k+n} u_{k,n}(r)^2\Theta_v(\d r) =\ff 1 {\ll_{k,n}} \int_k^{k+n} r \underline{a}(r) |u_{k,n}'(r)|^2 \Theta_v(\d r).$$
Letting $F_{k,n}(\eta)= u_{k,n}((\eta(E)\lor k)\land (k+n))$, for large enough $k>0$ such that $\GA_\eta=a(\rr(\eta))\1$ and $V(\eta)=v(\rr(\eta))$ for $\eta(E)\ge k$, the above formula implies
$$\BG^V(F_{k,n}^2)- \BG^V(F_{k,n})^2 \ge \BG^V(F^2_{k,n} \cdot 1_{\{k\le \rr\le k+n\}}) =\ff 1 {\ll_{k,n}} \EE_{\GA,V}(F_{n,k}, F_{n,k}),\ \ n\ge 1.$$
Obviously, due to \eqref{**} this implies $\gap(\L_{\GA,V})=0.$

  \end{proof}

  \subsection{Proof of Theorem \ref{NT}}

  Let $h\in \scr H$, i.e. $h\in C^1([0,\infty))$ with $h(r),h'(r)>0$ for $r>0$ such that \eqref{0*1} and \eqref{0*2} hold. By \eqref{PCW} and noting that $V_h=V-h(\rr)+c(h)$ where  $c(h)$ is a constant such that $\BG^{V_h}$ is a probability measure, for any $F\in\F C_0^\infty$ we have
  \beg{align*} &\BG^V(F^2\cdot 1_{\BB_N}) -\ff{\BG^V(F\cdot 1_{\BB_N})^2}{\BG^V(\BB_N)} =\inf_{c\in\R,|c|\le\|F\|_\infty} \BG^V(|F-c|^2\cdot 1_{\BB_N})\\
  &\le \e^{h(N)-c(h)} \inf_{c\in\R,|c|\le\|F\|_\infty} \BG^{V_h}(|F-c|^2\cdot 1_{\BB_N})\le \e^{h(N)-c(h)} \inf_{c\in\R,|c|\le\|F\|_\infty} \BG^{V_h}(|F-c|^2)\\
  &=  \e^{h(N)-c(h)} \big[\BG^{V_h}(F^2)- \BG^{V_h}(F)^2\big]
  \le C(h) \e^{h(N)-c(h)}\BG^{V_h} (\GG_\GA(F,F))\\
  &\le  C(h) \e^{h(N)}\BG^{V} (\GG_\GA(F,F))= C(h)\e^{h(N)} \EE_{\GA,V}(F,F).\end{align*}
 This implies
 $$\BG^V(F^2\cdot 1_{\BB_N})\le C(h)\e^{h(N)} \EE_{\GA,V}(F,F) + \ff{\BG^V(F\cdot 1_{\BB_N})^2}{\BG^V(\BB_N)},\ \ F\in\D(\EE_{\GA,V}).$$
  Then for any $F\in \D(\EE_{\GA,V})$ with $\BG^V(F)=0$, we have $\BG^V(F\cdot 1_{\BB_N})^2= \BG^V(F\cdot 1_{\{\rr>N\}})^2$ and
  \beg{align*}& \BG^V(F^2)\le \BG^V(F^2\cdot 1_{\BB_N}) +\BG^V(F^2\cdot 1_{\{\rr>N\}})\\
  &\le C(h)\e^{h(N)} \EE_{\GA,V}(F,F) + \ff{\BG^V(F\cdot 1_{\BB_N})^2}{\BG^V(\BB_N)}+\BG^V(F^2\cdot 1_{\{\rr>N\}})\\
  &\le  C(h)\e^{h(N)} \EE_{\GA,V}(F,F) + \Big(\ff{\BG^V(\rr>N)^2}{\BG^V(\BB_N)}+\BG^V(\rr>N)\Big) \|F\|_\infty^2 \\
  &\le  C(h)\e^{h(N)} \EE_{\GA,V}(F,F) +  \ff{ \BG^V(\rr>N)}{\BG^V(\BB_N)} \|F\|_\infty^2.\end{align*}
  So, for any $r>0$ and $N>0$ such that $\ff{\BG^V(\rr>N)}{\BG^V(\BB_N)}  \le r$, equivalently $\BG^V(\rr>N)\le \ff r {1+r}$, we have
  $$\BG^V(F^2)\le C(h)\e^{h(N)} \EE_{\GA,V}(F,F) + r\|F\|_\infty^2,\ \ F\in \D(\EE_{\GA,V}), \BG^V(F)=0, h\in \scr H.$$
  Therefore, the weak Poincar\'e inequality \eqref{WPC} holds for
  $$\aa(r):= \inf\Big\{C(h)\e^{h(N)}: \ h\in \scr H,  \BG^V(\rr>N) \le   \ff r{1+r}  \Big\},\ \ r>0.$$

\section{ Extensions to the space of finite signed measures}

Consider the space of finite signed measures
$$\MM_\s:= \{\eta-\eta':\ \eta,\eta'\in \MM\}$$ equipped with the  topology induced by the map
$$\eta\mapsto (\eta^+,\eta^-)\in\MM\times\MM,$$
where $\eta^+$ and $\eta^-$ are the positive and negative parts of $\eta$  in the Hahn decomposition respectively, and $\MM\times\MM$ is equipped with the weak topology. So, under this topology $\MM_\s$ is   a Polish space. Note that this topology maybe different from the weak topology,
 i.e. $\eta_n\to \eta$ if $\eta_n(f):=\int_E f\d\eta_n\to \eta(f)$ holds for any $f\in C_b(E)$, since the latter on $\MM_\s$ might be not metrizable, see \cite{Var}.

 To extend the Dirichlet form $(\EE_{\GA,V},\D(\EE_{\GA,V}))$  from $L^2(\BG^V)$ to  $L^2(\BG^V_\s)$ for a  probability measure $\BG^V_\s$ with a potential $V$ on $\MM_\s$,
 we   introduce below the measure $\BG^V_\s$,   the extrinsic derivative and   the operator $\GA$ respectively.

 In \cite{TVY}, an  analogue to the Lebesgue measure  was introduced on $\MM_\s$ by using the  convolution of two weighted Gamma distributions.   In the same spirit, we extend the measure $\BG$  to $\BG_\s$ on $\MM_\s$ as follows:
\beq\label{MC'} \int_{\MM_\s} f( \eta) \BG_\s(\d \eta)= \int_{\MM\times\MM} f(\eta^+-\eta^-)  \BG(\d\eta^+)\BG(\d\eta^-),\ \ f\in \B_b(\MM_\s).\end{equation}   Let $\tau(\eta)= \{x\in E: \eta(\{x\})\ne 0\}$. To ensure that $\tau(\eta^+)$ and $\tau(\eta^-)$ are disjoint such that $\eta=\eta^+-\eta^-$ is the Hahn decomposition of $\eta$, we will assume that $\nu$ is atomless. In this case,  $\tau(\eta^+)\cap\tau(\eta^-)=\emptyset $  for $\BG\times\BG$-a.e. $(\eta^+,\eta^-)$.

 Next, we define the extrinsic derivative operator  $(\nn^{ext}, \D(\nn^{ext}))$ as  in Definition 1.1 for $\MM_\s$ replacing $\MM$:
 \beq\label{*PW} \nn^{ext} F(\eta)(x) = \lim_{0\ne s\to  0} \ff{F(\eta+s\dd_x)-F(\eta)}s,\ \ \eta\in\MM_\s.\end{equation}
Let   $ \F_\s C_0^\infty$ be the class of cylindrical functions of type
\beq\label{*PY0-1}   F(\eta):=f(\eta^+(A_1),\cdots, \eta^+(A_n),\eta^-(A_1),\cdots, \eta^-(A_n) ),\ \ n\ge 1, f\in C_0^\infty(\R^{2n}),\end{equation}  where $  \{A_i\}_{1\le i\le n}$ is  a measurable partition  of $E$, and $\eta=\eta^+-\eta^-$ is the Hahn decomposition.  Let
\beq\label{*PY0} A_\eta:= \{x\in E:\ \eta(\{x\})<0\},\ \  \eta\in\MM_\s.\end{equation}
It is easy to see that such a function $F$ is extrinsically differentiable with
\beq\label{*PY}\beg{split}    \nn^{ext}  F(\eta)= \   &\sum_{i=1}^{2n} (1-2\cdot 1_{\{i>n\}}) (\pp_if)(\eta^+(A_1),\cdots, \eta^+(A_n), \eta^-(A_1),\cdots, \eta^-(A_n)) 1_{A_\eta^{i,n}},\end{split}\end{equation} where
$$A_{\eta}^{i,n}:= \beg{cases} A_i\cap A_\eta^c,\ &\text{if}\ i\le n,\\
A_i\cap A_\eta,\ &\text{if} \ i>n.\end{cases}$$
Since   for any $\eta\in\MM_\s$,   $A_{\eta+\vv\dd_x}=A_\eta$ holds for small $\vv>0$ and all $x\in E$, $\nn^{ext} F(\eta)(x)$ is again extrinsically differentiable in $\eta$ with
\beq\label{*PY'}\beg{split}   &\nn^{ext}[\nn^{ext}  F(\eta)(x)] (y)
  =  \sum_{i,j=1}^{2n}\Big[(1-2\cdot 1_{\{i>n\}})(1-2\cdot 1_{\{j>n\}}) \\
   &\qquad\qquad\times (\pp_i\pp_jf)(\eta^+(A_1),\cdots, \eta^+(A_n), \eta^-(A_1),\cdots, \eta^-(A_n)) 1_{A_\eta^{i,n}}(x) 1_{A_\eta^{j,n}}(y)\Big].\end{split}\end{equation}
Finally, For any $\eta\in \MM_\s$, let   $\GA_{\eta}$ be a positive definite bounded linear operator on $L^2(|\eta|)$, where $|\eta|:=\eta^++\eta^-$ is the total variation of $\eta$.
Consider the pre-Dirichlet form
\beq\label{PDs}  \EE_{\GA,V}^\s(F,G):= \int_{\MM_\s} \<\GA_{\eta} \nn^{ext}F(\eta), \nn^{ext}G(\eta)\>_{L^2(|\eta|)} \,\d\BG^V_\s,\ \ F,G\in \F C_0^\infty.\end{equation}
To ensure the closability of this bilinear form,     we assume

\beg{enumerate} \item[ {\bf (H')}]  $\nu$ is atomless, $V\in\D(\nn^{ext})$ such that $\BG^V_\s$ is a probability measure. Moreover, for any
  $A\in\B(E)$ and $x\in E$,  $\GA_\eta 1_{A\cap A_\eta^c}(x)$   and   $\GA_\eta 1_{A\cap A_\eta}(x)$  are  extrinsically   differentiable in  $\eta$  with
\beg{align*}
  &\int_{\MM_\s} \Big[  |\eta|\big(|\nn^{ext}[\GA_\eta 1_{A\cap A_\eta^c}]| +|\nn^{ext}[\GA_\eta 1_{A\cap A_\eta}]|\big)\\
  &\qquad\qquad + |\eta|\big((|\GA_\eta 1_{A\cap A_\eta^c}|+|\GA_\eta 1_{A\cap A_\eta}|)|\nn^{ext}V(\eta)|  \big) \Big]\BG^V_\s(\d\eta)<\infty. \end{align*}
  \end{enumerate} Obviously, this assumption is satisfied  if $\GA_\eta=F(\eta)\1$ for some positive bounded  extrinsically differentiable function $F$ such that  $\BG^V_\s$ is a probability measure with
  $$\int_{\MM_\s} |\eta|(|\nn^{ext}F(\eta)|+|\nn^{ext}V(\eta)|)\BG^V_\s(\d\eta)<\infty.$$

\subsection{Integration by parts formula}

 \beg{thm}\label{T4.1} Assume {\bf (H')}.  Then
 \beq\label{DRC} \EE_{\GA,V}^\s(F,G)=-\int_{\MM_\s}( G\L_{\GA,V}^\s F) \, \d\BG^V_\s,\ \ F,G\in \F_{\rm s} C_0^\infty \end{equation} holds for
\beg{align*} & \L_{\GA,V}^\s F(\eta)   :=  \int_E \Big(\big[\nn^{ext}V(\eta)(x)\big]\GA_\eta[\nn^{ext}F(\eta)](x) + \nn^{ext}\big(\GA_\eta[\nn^{ext}F(\eta)](x)\big)(x) \Big)\,|\eta|(\d x)\\
  &\qquad \qquad \qquad -\int_E \GA_\eta[\nn^{ext} F(\eta)](x) \eta(\d x).\end{align*}
 Consequently, $(\EE_{\GA,V}^\s, \F_{\rm s} C_0^\infty)$ is closable in $L^2(\BG^V_\s)$ and its closure $(\EE_{\GA,V}^\s, \D(\EE_{\GA,V}^\s))$ is a symmetric Dirichlet form with $1\in\D(\EE_{\GA,V}^\s)$ and $\EE_{\GA,V}^\s(1,1)=0.$
 \end{thm}

To prove this result, we   introduce the divergence operator associated with $\nn^{ext}$.

\beg{defn}  A  measurable  function $\phi $ on $ \MM_\s\times E$ is said in the domain $\D(\ddiv^{ext}_\s)$,
if for any $x\in E$ we have $\phi(\cdot,x)\in\D(\nn^{ext})$ and
\beq\label{BBT2}  \int_{\MM_\s} \bigg( \int_E\big(|\nn^{ext} \phi(\eta,x)(x)| +|\phi(\eta,x)| \big)|\eta|(\d x)\bigg)\, \BG_\s(\d\eta)<\infty.\end{equation}
In this case,  the divergence operator is given  by
\beq\label{DIVs} \ddiv^{ext}_\s(\phi)(\eta):= \int_E  \phi(\eta,x)\eta(\d x)- \int_E \nn^{ext}\phi(\eta,x)(x)\,|\eta|(\d x),\ \ \eta\in\MM_\s.\end{equation}
\end{defn}
We have the following integration by parts formula for the directional derivative
 $$\nn_\phi^{ext}F(\eta):= \int_E\big[ \phi(\eta, x) \nn^{ext}F(\eta)(x)\big]|\eta|(\d x),\ \ \phi\in \D(\ddiv^{ext}_\s), F\in\D(\nn^{ext}).$$

 \beg{lem}\label{L4.2}   Let  $\phi\in \D(\ddiv^{ext}_\s)$. Then    $$\int_{\MM_\s} (\nn^{ext}_\phi F) \, \d\BG_\s  =  \int_{\MM_\s} [F\ddiv^{ext}_\s(\phi)] \,\d\BG_\s, \ \ F\in\F_\s C_0^\infty.$$
  \end{lem}
    \beg{proof}  By a simple approximation argument, we may and do assume that $\phi$ is bounded so that $(\BG_\s\times\nu)(|\phi|)<\infty$.
    For $F\in \F_{\rm s} C_0^\infty$, \eqref{*PW}  implies
  \beq\label{ABC} \nn^{ext}F(\eta)(x)= \nn^{ext}F(\cdot-\eta^-)(\eta^+)(x)=- \nn^{ext} F(\eta^+-\cdot)(\eta^-)(x),\ \ F\in\D(\nn^{ext}).\end{equation}
 Next,   for any $\eta'\in\MM$, let
  $$\phi_{+,\eta'}(\eta,x) := \phi(\eta'-\eta,x),\ \ \phi_{-,\eta'}(\eta,x) := \phi(\eta-\eta',x),\  \ (\eta,x)\in \MM\times E,$$
By \eqref{DIV} and \eqref{ABC}  we obtain
\beg{align*}& \ddiv^{ext}(\phi_{-,\eta^-})(\eta^+)- \ddiv^{ext}(\phi_{+,\eta^+})(\eta^-)\\
&= \int_E \big[\phi(\eta^+-\eta^-,x)-\nn^{ext}\phi(\cdot-\eta^-,x)(\eta^+)(x)\big]\eta^+(\d x)-\nu(\phi(\eta,\cdot))\\
&\qquad - \int_E \big[\phi(\eta^+-\eta^-,x)+ \nn^{ext}\phi(\eta^+-\cdot,x)(\eta^-)(x)\big]\eta^-(\d x)-\nu(\phi(\eta,\cdot)) \\
&= \int_E\phi(\eta^+-\eta^-,x)(\eta^+-\eta^-)(\d x) - \int_E \big[\nn^{ext}(\cdot,x)(\eta^+-\eta^-)(x)\big] (\eta^++\eta^-)(\d x)\\
&=\ddiv^{ext}_\s(\phi)(\eta),\ \ \eta=\eta^+-\eta^-\ \text{with}\ \tau(\eta^+)\cap\tau(\eta^-)=\emptyset.\end{align*}
  Combining this with  Lemma \ref{L2.2}, \eqref{MC'} and \eqref{ABC},  we obtain
\beg{align*} &\int_{\MM_\s} (\nn^{ext}_\phi F) \, \d\BG_\s= \int_{\MM\times\MM}\BG(\d\eta^+) \BG(\d\eta^-)\int_E \big[\phi(\eta^+-\eta^-,x) \nn^{ext}F(\eta^+-\eta^-)(x)\big] (\eta^++\eta^-)(\d x) \\
&= \int_{\MM} \BG(\eta^-) \int_{\MM \times E}  \big[\phi(\eta^+-\eta^-,x) \nn^{ext}F(\cdot-\eta^-)(\eta^+)(x)\big] \eta^+ (\d x) \\
&\qquad -\int_{\MM} \BG(\eta^+) \int_{\MM \times E}  \big[\phi(\eta^+-\eta^-,x) \nn^{ext}F(\eta^+-\cdot)(\eta^-)(x) \big] \eta^- (\d x)\\
&=  \int_{\MM\times\MM} F(\eta^+-\eta^-)\big[ \ddiv^{ext} (\phi_{-,\eta^-})(\eta^+) -\ddiv^{ext} (\phi_{+,\eta^+})(\eta^-) \big]  \BG(\d\eta^+)\BG(\eta^-) \\
&= \int_{\MM_\s} F(\eta) \ddiv^{ext}_\s(\phi)(\eta) \,\BG_\s(\d\eta).\end{align*}

   \end{proof}

\beg{proof}[Proof of Theorem \ref{T4.1}]   Let $F\in \F_\s C_0^\infty$ be given in \eqref{*PY0-1}, and let
\beg{align*} & \phi(\eta, x):=\e^{V(\eta)}\GA_\eta [\nn^{ext} F(\eta)] (x)\\
&= \e^{V(\eta)} \sum_{i=1}^{2n } (1-2\cdot 1_{\{i>n\}}) (\pp_i f)(\eta(A_1,\cdots, \eta(A_n)) \GA_\eta 1_{A_i^\eta}(x),\ \ (\eta,x)\in \MM_\s\times E.\end{align*}Then   {\bf (H')} and \eqref{*PY} imply $\phi\in\D(\ddiv_\s^{ext})$.
By the definition of $\EE_{\GA,V}^\s$ and Lemma \ref{L4.2}, for any $G\in\F_\s C_0^\infty$ we have
\beg{align*} &\EE_{\GA,V}^\s(F,G)=    \int_{\MM_\s} \big\<\GA_\eta \nn^{ext} F(\eta), \nn^{ext} G(\eta) \big\>_{L^2(|\eta|)}  \,\BG_\s^V(\d\eta)\\
  &= \int_{\MM_\s}  \big\<\phi(\eta,\cdot), \nn^{ext}G(\eta)\big\>_{L^2(|\eta|)} \, \BG_\s(\d\eta)  =  \int_{\MM_\s}G(\eta) \ddiv^{ext}_\s(\phi)\,\BG_\s(\d\eta).\end{align*}
This together with  \eqref{DIVs} implies \eqref{DRC}   for
  \beg{align*} &\L_{\GA,V}^\s F(\eta):= -\e^{-V(\eta)}\ddiv^{ext}_\s(\phi)= - \e^{-V(\eta)} \ddiv^{ext}_\s\big(\e^{V(\eta)}\GA_\eta[\nn^{ext}F(\eta)](\cdot)\big)\\
  &=  \int_E \Big(\big[\nn^{ext}V(\eta)(x)\big]\GA_\eta[\nn^{ext}F(\eta)](x) + \nn^{ext}\big(\GA_\eta[\nn^{ext}F(\eta)](x)\big)(x) \Big)\,|\eta|(\d x)\\
  &\qquad -\int_E \GA_\eta[\nn^{ext} F(\eta)](x) \eta(\d x).\end{align*}

  Next, to prove that $1\in\D(\EE_{\GA,V}^\s)$ with $\EE_{\GA,V}^\s(1,1)=0$, we take   $\{f_n\}_{n\ge 1}\subset C_0^\infty(\R)$ such that $f_n(s)=1$ for $|s|\le n$, $0\le f_n\le 1$ and $ \|f_n'\|_\infty\le 1$.
Let $F_n(\eta):= f_n(\eta(E)), n\ge 1$. Then $F_n\in\F C_0^\infty$. By    {\bf (H')} we have     $\BG^V_\s(|F_n-1|^2)\to 0$ as $n\to\infty$, and
$$\limsup_{n\to\infty} \EE_{\GA,V}^\s(F_n,F_n) = \limsup_{n\to\infty} \int_{\{|\eta(E)|>n\}} \|\GA_\eta 1\|_{L^1(|\eta|)} \, \BG^V_\s(\d\eta) =0.$$
Therefore, $1\in \D(\EE_{\GA,V}^2)$ and $\EE_{\GA,V}^\s(1,1)=0.$
 \end{proof}

 \subsection{Functional inequalities for $\EE_{\1,0}^\s$}
 For any $N>0$, let $\tt \BB_N^\s=\{\eta\in \MM_\s: \eta^+(E)\lor \eta^-(E)\le N\}$.

 \beg{thm}\label{T5.2} Let $\GA=\1$ and $V=0$.
 \beg{enumerate} \item[$(1)$]  $\gap(\L_{\GA,V}^\s)= 1$, i.e. the following Poincar\'e inequality
\beq\label{P00} \BG_\s(F^2)\le    \EE_{\1,0}^\s(F,F)+ \BG_\s(F)^2,\ \ F\in \D(\EE_{\1,0}^\s )\end{equation}  holds, and the constant $1$ in front of $\EE_{\1,0}^\s(F,F)$ is optimal.
 \item[$(2)$] If  ${\rm supp}\,\nu$ is infinite, then $\EE_{\1,0}^\s$ does not satisfy the super Poincar\'e inequality. On the other hand,    there exists a constant $c_0>0$ such that when ${\rm supp}\,\nu$ is a finite set,   the log-Sobolev inequality
\beq\label{LSI0*} \BG_\s(F^2\log F^2)\le \ff{c_0}{ 1\land \dd}\,\EE_{\1,0}^\s(F,F),\ \ F\in \D(\EE_{\1,0}), \BG(F^2)=1\end{equation}
holds, where $\dd:= \min\{\nu(\{x\}):\ x\in  {\rm supp}\,\nu\}.$
\item[$(3)$]  For any $N>0$ and $F\in \F C_0^\infty$ with $\BG_\s(1_{\tt \BB_N^\s} F)=0,$
 \beq\label{LSP*} \BG_\s(1_{\tt \BB_N^\s} F^2 ) \le  \Big(2\lor\ff{N^2}{2\nu(E)}\Big)  \BG_\s\big(1_{\tt B_N^\s} \|\nn^{ext} F\|_{L^2(|\eta|)}^2 \big).\end{equation}\end{enumerate}
  \end{thm}
 \beg{proof}   By taking $F(\eta)$ depending only on $\eta^+$, it is easy to see that a Poincar\'e inequality for $\EE_{\1,0}^\s$ implies the same inequality for $\EE_{\1,0}$. So, the optimality of \eqref{P00}, and the invalidity of the super Poincar\'e inequality when supp\,$\nu$ is infinite, follow from Theorem \ref{T1.2}. It remains   to prove the inequalities
 \eqref{P00}, \eqref{LSI0*} and \eqref{LSP*}.    According to the additivity property of the Poincar\'e  and log-Sobolev inequalities,
    these inequalities follow  from the corresponding ones of $\EE_{\1,0}$.   For simplicity, below we only prove the first inequality.

 Let $\F\in \F_\s C_0^\infty$. By Theorem \ref{T1.2}, \eqref{MC'},  \eqref{PDs} for $\GA=\1$ and $V=0$, and using \eqref{ABC}, we obtain
 \beg{align*} &\BG_\s( F^2) = \int_\MM \BG(\d\eta^-) \int_\MM F(\eta^+-\eta^-)^2 \BG(\d\eta^+) \\
 &\le       \int_{\MM \times\MM} \|\nn^{ext} F(\cdot-\eta^-)(\eta^+)\|^2_{L^2(\eta^+)}  \BG(\d\eta^+)\BG(\d\eta^-)+\int_\MM \bigg(\int_\MM F(\eta^+-\eta^-)  \BG(\d\eta^+)\bigg)^2  \BG(\d\eta^-)  \\
 &\le   \int_{\MM_\s} \|\nn^{ext}F(\eta)\|^2_{L^2(\eta^+)}\BG_\s(\d\eta)+ \bigg( \int_{\MM \times\MM}   F(\eta^+-\eta^-)  \BG(\d\eta^+)\BG(\d\eta^-)\bigg)^2  \\
 &\qquad +  \int_\MM \bigg\|\nn^{ext} \bigg[\int_\MM F(\eta^+-\cdot)\BG(\d\eta^+)\bigg](\eta^-)\bigg\|_{L^2(\eta^-)}^2 \BG(\d\eta^-).\end{align*}
 By the Jensen inequality, we have
$$\bigg\|\nn^{ext} \bigg[\int_\MM F(\eta^+-\cdot)\BG(\d\eta^+)\bigg](\eta^-)\bigg\|_{L^2(\eta^-)}^2 \le \int_{\MM} \|\nn^{ext} F(\eta^+-\cdot)(\eta^-)\|_{L^2(\eta^-)}^2 \BG(\d\eta^+).$$
Therefore,
$$\BG_\s( F^2)  \le \BG_\s(F)^2 +  \int_{\MM_\s} \|\nn^{ext} F(\eta)\|^2_{L^2(\eta^++\eta^-)} \BG_\s(\d\eta)= \BG_\s(F)^2+  \EE_{\1,0}^\s(F,F).$$

   \end{proof}

  \subsection{Functional inequalities for $\EE_{\GA,V}^\s$}

According to the proof of Theorem \ref{TSP} and the local Poincar\'e inequality \eqref{LSP*}, it seems that we should take
$$\tt\rr_\s(\eta):=2\sqrt{\eta^+(E)\lor \eta^-(E)},\  \ \eta\in\MM_\s$$ to replace the function $\rr$ on $M$.  But by \eqref{*PY}
we have $$\nn^{ext} \tt\rr_\s (\eta)(x)= \ff 2 {\tt\rr_\s(\eta)}\Big( 1_{\{\eta(E)\ge 0\}}1_{A_\eta^c}(x)- 1_{\{\eta(E)< 0\}}1_{A_\eta}(x)\Big),$$ which is however not extrinsically differentiable in $\eta$, so that  $\L_{\GA,V}\tt\rr_\s$ is not well defined as required. To avoid this problem, below we will use both $\tt\rr_\s$ and
\beq\label{BBG} \rr_\s(\eta):= 2\sqrt{|\eta|(E)},\ \ \eta\in\MM_\s,\end{equation}
which satisfies $\|\nn^{ext} \rr_\s(\eta)\|_{L^2(|\eta|)}=1$ according to the following lemma.

\beg{lem}\label{L5.1} Let $\rr$ be defined in $\eqref{BBG}$ and let $\s(\eta,\cdot):= 1-2\cdot 1_{A_\eta}$ for $A_\eta$ in $\eqref{*PY0}$.  Then
\beg{align*}  \nn^{ext}\rr_\s(\eta)= \ff {2\s(\eta,\cdot)} {\rr_\s(\eta)},\ \
\nn^{ext}\s(\eta,x)(y)=0,\ \ \eta\in\MM_\s, x,y\in E.\end{align*}
Consequently,
 if $\GA_\eta \s(\eta,\cdot)$ is extrinsically differentiable in $\eta\in\MM_\s$ with
  \beq\label{LOT} \sup_{|\eta|(E)\le r}  |\eta|\big( |\GA_\eta \s(\eta,\cdot)|+ |\nn^{ext} [\GA_\eta\s(\eta,\cdot)](\cdot)| \big) <\infty,\ \  r\in (0,\infty),\end{equation}
then
\beq\label{LRR}  \L_{\GA,V}^\s \rr_\s(\eta)
   =  \ff 2 {\rr_\s(\eta)} \Big[|\eta|\Big([\nn^{ext}V(\eta)] \GA_\eta\s(\eta,\cdot)\big) +  \nn^{ext}\big[\GA_\eta \s(\eta,\cdot)\big](\cdot)\Big)
  - \eta\big(\GA_\eta \s(\eta,\cdot)\big)\Big].
  \end{equation}
\end{lem} This  lemma can be proved by simple calculations using \eqref{*PW} and the definition of $\L_{\GA,V}$  in Theorem \ref{T4.1}, so we omit the details.

By Lemma \ref{L5.1},   we have $$\GG_\1^\s(\rr_\s,\rr_\s):= |\eta|(|\nn^{ext}\rr_\s|^2)=1,\ \ \nn^{ext}[\nn^{ext}\rr_\s(\eta)(x)](x)= -\ff{4}{\rr_\s(\eta)^3}.$$ These coincide with the
corresponding   properties of $\rr$ on $\MM$.

Similarly to \eqref{*TP0} and \eqref{*TP1},  let
\beq\label{*FY} \beg{split} &\xi_\s(r)= \inf_{\rr_\s(\eta)=r} \L_{\GA,V}^\s\rr_\s(\eta), \ \  \underline{a}_\s(r)= \inf_{\rr_\s(\eta)=r}\inf_{\|\phi\|_{L^2(|\eta|)=1} } \<\GA_\eta\phi,\phi\>_{L^2(|\eta|)},\\
&\bar a_\s(r)= \sup_{\rr_\s(\eta)=r}\sup_{\|\phi\|_{L^2(|\eta|)=1} } \<\GA_\eta\phi,\phi\>_{L^2(|\eta|)},\ \ \ r>0,\\
&\si_{k,\s}=   \sup_{t\ge k} \int_t^\infty   \e^{\int_k^r \ff{ \xi_\s(s)}{\underline{a}_\s(s)} \d s} \d r\int_k^t  \ff 1 {  \underline{a}_\s(r)} \e^{-\int_k^r \ff{ \xi_\s (s)}{\underline{a}_\s(s)} \d s} \d r,\ \ k>0.
 \end{split}\end{equation}
Assume that
\beq\label{00*} \psi(t):= \int_0^t [\bar a_\s(r)]^{-\ff 1 2}\,\d  r\uparrow\infty\ \text{as}\ t\uparrow\infty.\end{equation}
As in the proof of Theorem \ref{TSP},   we may use $ \si_{k,\s}$    to estimate $\BG_\s^V(F_N^2)$ for
$$F_N:= \big[(\psi(\rr_\s)-\psi(N))^+\land 1\big]\cdot  F,\ \ N>0, F\in \F_\s C_0^\infty.$$
More precisely, as in \eqref{PPK0} and \eqref{PPK} we conclude that for any $k>0$ there exists $N\in [k, \psi^{-1}(\psi(k)+32\si_{k,\s})]$   such that
\beq\label{PPK2} \beg{split} \int_{\MM_\s} F_N^2 \d \BG^V_\s&\le \ff 2 { \ll_k} \EE_{\GA,V}^\s(F,F) + \ff 1 4 \int_{\MM} F^2 \d\BG^V_\s\\
&\le 8\si_{k,\s}\, \EE_{\GA,V}^\s(F,F) + \ff 1 4 \int_{\MM} F^2 \d\BG^V_\s.\end{split} \end{equation}
On the other hand, we estimate $\BG^V_\s(F^2\cdot 1_{\{\rr_\s\le N\}})$ by using the local Poincar\'e inequality \eqref{LSP*}.
Since   the bounded set in \eqref{LSP*} is   $\tt\BB_N^\s:=\{\tt\rr_\s\le N\}$ rather than $\BB_N^\s:=\{\rr_\s\le N\}$,
    we change the definition of $\Phi(N)$ into
\beg{align*} \Phi_{\s}(N):=\Big(2\lor
\ff{N^2}{2\nu(E)}\Big) \exp\bigg[\sup_{\tt\rr_\s\le N}V-\inf_{\tt\rr_\s\le N}V\bigg] \sup_{\tt\rr_\s(\eta)\le N}\sup_{\|\phi\|_{L^2(|\eta|)}=1} \ff 1 {\<\GA_\eta \phi,\phi\>_{L^2(|\eta|)}},\ \ N>0.\end{align*}
Noting that   $ 1_{\{\rr_\s\le N\}}\le  1_{\{\tt\rr_\s\le N\}}$,   we may apply
Theorem \ref{T5.2} to bound  $\BG^V_\s(F^2\cdot 1_{\{\rr_\s\le N\}})$.  For instance, corresponding to \eqref{*A5} we have
$$\BG^V_\s(F^2\cdot 1_{\{\rr_\s\le N\}})\le \BG^V_\s(F^2\cdot 1_{\{\tt\rr_\s\le N\}})\le  \BG_\s^V(1_{\{\tt\rr_\s\le N\}}F)^2 +   \Phi_{\s}(N)\EE_{\GA,V}^\s( F,F).$$
Combining this with \eqref{PPK2} we may extend assertions of Theorem \ref{TSP} to the present setting as follows, where when ${\rm supp}\,\nu$  is infinite the super Poincar\'e can be disproved as in the proof of Theorem \ref{TSP}(2) by taking $F_n(\eta)= (1-\eta^+(E))^+\ff{\eta^+(A_n)}{\eta^+(E)}$ for $0<\nu(A_n)\downarrow 0$. Moreover, one may also extend Corollaries $\ref{C1.2}$-$\ref{C2.4}$ and Theorem \ref{NT}.
   We omit the details to save space.

\beg{thm}\label{TSP*} In addition to {\bf(H')}, assume that $\GA_\eta\s(\eta,\cdot)$ is extrinsically differentiable in $\eta$ such that $\eqref{LOT}$ holds. Moreover, assume that $\underline{a}_\s$ and $\bar a_\s$ in $\eqref{*FY}$ are such that $\underline{a}_\s^{-1}(r)$ is locally bounded in $r\ge 0$ and $\eqref{00*}$ holds.
   \beg{enumerate}
\item[$(1)$] If   $\lim_{k\to\infty}   \si_{k,\s}<\infty$, then
 $$\gap(\L_{\GA,V}^\s )\ge \sup\bigg\{\ff 1 {2 \Phi_\s\big(\psi^{-1}(\psi(k)+ 32  \si_{k,\s} + 1 )\big) + 32   \si_{k,\s}}:\ k>0\bigg\}>0.$$
 \item[$(2)$]  If ${\rm supp}\,\nu$ contains  infinitely many points, then $\EE_{\GA,V}^\s$ does not satisfy the super Poincar\'e inequality.
  \item[$(3)$] The weak Poincar\'e inequality $\eqref{WPC}$  holds for $(\EE_{\GA,V}^\s,\BG_\s^V)$ replacing $(\EE_{\GA,V},\BG^V)$ and
 $$\aa(r):= \inf\Big\{2 \Phi_\s(N):\  \BG^V_\s(\rr>N) \le \ff r{1+r}\Big\},\ \ r>0.$$ \end{enumerate}
\end{thm}

 \paragraph{Acknowledgement.} The author  would like to thank Professor Shui Feng for valuable conversations,
as well as  the referee for very careful reading, helpful comments and corrections.

\end{document}